\documentclass[11pt]{article}
\usepackage{amsmath}
\usepackage{dsfont}
\usepackage{mathrsfs}
\usepackage{amsmath,amssymb}
\usepackage{amsfonts}
\usepackage{hyperref}
\usepackage{amsthm}
\usepackage{graphicx}
\usepackage{subfigure}
\usepackage{xcolor}
\renewcommand{\qed}{\hfill\small{$\blacksquare$}\normalsize}

\hfuzz=\maxdimen
\theoremstyle{definition}
\newtheorem{lemma}{Lemma}[section]
\newtheorem{definition}[lemma]{Definition}
\newtheorem{proposition}[lemma]{Proposition}
\newtheorem{theorem}[lemma]{Theorem}
\newtheorem{corollary}[lemma]{Corollary}
\newtheorem{example}{Example}
\newtheorem{remark}{Remark}

\numberwithin{equation}{section}
\renewcommand{\proof}{\textbf{Proof. }}
\renewcommand{\qed}{\hfill\small{$\blacksquare$}\normalsize}

\DeclareFixedFont{\Acknowledgment}{OT1}{cmr}{bx}{n}{14pt}
\begin{document}

\title{\bf On a combinatorial curvature for surfaces with inversive distance circle packing metrics}
\author{Huabin Ge, Xu Xu}
\maketitle

\begin{abstract}
In this paper, we introduce a new combinatorial curvature on triangulated surfaces with inversive
distance circle packing metrics. Then we prove that this combinatorial curvature has
global rigidity. To study the Yamabe problem of the new curvature,
we introduce a combinatorial Ricci flow, along which the curvature evolves almost in the same way as
that of scalar curvature along the surface Ricci flow obtained by Hamilton \cite{Ham1}.
Then we study the long time behavior of the combinatorial Ricci flow
and obtain that the existence of a constant curvature metric is
equivalent to the convergence of the flow on triangulated surfaces with nonpositive Euler number.
We further generalize the combinatorial curvature to $\alpha$-curvature and prove that it is also globally
rigid, which is in fact a generalized Bower-Stephenson conjecture \cite{BS}. We also use the combinatorial Ricci
flow to study the corresponding $\alpha$-Yamabe problem.
\end{abstract}


\section{Introduction}\label{Introduction}

This is a continuation of our work on combinatorial curvature in \cite{GX4}.
This paper generalizes our results in \cite{GX4} to triangulated surfaces
with inversive distance circle packing metrics.
Circle packing is a powerful tool in the study of differential geometry and geometric topology and there are
lots of research on this topic.
In his work on constructing hyperbolic structure on 3-manifolds, Thurston (\cite{T1}, Chapter 13) introduced
the notion of Euclidean and hyperbolic circle packing metrics on triangulated surfaces
with prescribed intersection angles.
The requirement of prescribed intersection angles corresponds to the fact that the intersection angle of
two circles is invariant under the M\"{o}bius transformations.
For triangulated surfaces with Thurston's circle packing metrics, there will be singularities at the vertices.
The classical combinatorial Gauss curvature $K_i$ is introduced to describe the singularity at the vertex $v_i$,
which is defined as the angle deficit at $v_i$.
Thurston's work generalized Andreev's work on circle packing metrics on a sphere \cite{An1,An2}. Andreev and
Thurston's work together gave a complete characterization of the space of the classical combinatorial Gauss curvature.
As a corollary, they got the combinatorial-topological obstacle for the existence of a constant
curvature circle packing metric, which could be written as combinatorial-topological inequalities.
Chow and Luo \cite{CL1} first introduced a combinatorial Ricci flow, a combinatorial analogue of the smooth
surface Ricci flow, for triangulated surfaces with Thurston's circle packing metrics and got the equivalence between the
existence of a constant curvature metric and the convergence of the combinatorial Ricci flow.
This work is the cornerstone of applications of combinatorial surface Ricci flow in engineering up to now, see for example \cite{ZG, ZZGLG} and
the references therein.
Luo \cite{L1} once introduced a combinatorial Yamabe flow on triangulated surfaces with piecewise linear metrics to
study the corresponding constant curvature problem.
The combinatorial surface Ricci flow in \cite{CL1} and the combinatorial Yamabe flow in \cite{L1} are recently written
in a unified form in \cite{ZGZLYG}.
The first author \cite{Ge-thesis,Ge} introduced a combinatorial Calabi flow on triangulated surfaces with Thurston's Euclidean circle packing metrics
and proved the equivalence between the existence of constant circle packing metric and the convergence of the combinatorial Calabi flow.
The authors \cite{GX3} further generalized the combinatorial Calabi
flow to hyperbolic circle packing metrics and got similar results.

However, there are some disadvantages for the classical discrete Gauss curvature as stated in \cite{GX4}.
The first is that the classical Euclidean discrete Gauss curvature is invariant under scaling,
i.e., $K_i(\lambda r)=K_i(r)$ for any positive constant $\lambda$.
The second is that the classical discrete Gauss curvature tends to zero, not the Gauss curvature of the smooth surface,
as triangulated surfaces approximate a smooth surface.
Motivated by the two disadvantages, the authors \cite{GX4} introduced a new combinatorial curvature defined as $R_i=\frac{K_i}{r_i^2}$
for triangulated surfaces with Thurston's Euclidean circle packing metrics. If we take $g_i=r_i^2$ as the analogue of the Riemannian metric,
then we have $R_i(\lambda g)=\lambda^{-1}R_i(g)$, which has the same form as that of the smooth Gauss curvature.
Furthermore, there are examples showing that this curvature actually approximates the smooth Gauss curvature on surfaces as
the triangulated surfaces approximate a smooth surface \cite{GX4}.
Then we introduce a combinatorial Ricci flow and a combinatorial Calabi flow to study the corresponding Yamabe problem and got
a complete characterization of  the existence of a constant $R$-curvature circle packing metric using these flows.
The results in \cite{GX4} generalized the previous work of the authors in \cite{GX2}.
The authors then further generalized the curvatures to triangulated surfaces with hyperbolic background geometry \cite{GX3}.
We also consider the $\alpha$-curvature and $\alpha$-flows on low dimensional triangulated manifolds in \cite{GX4, GX5}.

Thurston introduced the intersection angle in the definition of circle packing metric as it is invariant under
the M\"{o}bius transformations. The notion of inversive distance of two circles in a M\"{o}bius plane was introduced
by H.S.M. Coxeter \cite{Co}, which generalizes the notion of the intersection angle of two circles.
Furthermore, the notion of inversive distance is invariant under the M\"{o}bius transformations, as it is actually
defined by the cross ratio \cite{BH}.
Bowers and Stephenson \cite{BS} introduced the inversive distance circle packing metric for triangulated surfaces, which
generalizes Thurston's circle packing metric defined using intersection angle.
The inversive distance circle packing metric has some applications in medical science, see for example \cite{HS}.
Bowers and Stephenson \cite{BS} once conjectured that the inversive distance circle packings are rigid.
This conjecture was proved to be locally right by Guo \cite{Guo} and then finally proved by Luo \cite{L2} in
Euclidean and hyperbolic background geometry for nonnegative inversive distance.
The local rigidity \cite{Guo} comes from the convexity of an energy function defined on the admissible space and
the global rigidity \cite{L2} comes from a convex extension of the energy function to the whole space, the idea of which
comes from Bobenko, Pinkall and Springborn \cite{BPS}.
Ma and Schlenker \cite{MS} gave a counterexample showing that there is no rigidity for the spherical background geometry
and John C. Bowers and Philip L. Bowers \cite{BB} recently presented a new construction
of their counterexample using only the inversive geometry of the 2-sphere.
The combinatorial Ricci flow for triangulated surfaces with inversive distance circle packing metrics
was introduced in  \cite{ZG, ZGZLYG}, using the classical combinatorial Gauss curvature.
Recently, the first author and Jiang \cite{GJ1, GJ2}studied the long time behavior of the flow using the extending method
and got some interesting results.

In this paper, we generalize our definition of combinatorial curvature in \cite{GX4} for Thurston's circle packing metrics
to inversive distance circle packing metrics. This curvature was also studied at the same time by the first author and Jiang \cite{GJ3}
using the flow introduced by the authors in \cite{GX5}.
Given a weighted triangulated surface $(M,\mathcal{T}, I)$ with inversive distance $I_{ij}\geq0$ for every edge $\{ij\}\in E$,
for any function $r: V\rightarrow [0, +\infty)$, the edge length $l_{ij}$ is defined to be
\begin{equation*}
l_{ij}=\sqrt{r_i^2+r_j^2+2r_ir_jI_{ij}}
\end{equation*}
for Euclidean background geometry and
\begin{equation*}
l_{ij}=\cosh^{-1}\left(\cosh r_i\cosh r_j+I_{ij}\sinh r_i\sinh r_j\right)
\end{equation*}
for hyperbolic background geometry. If the lengths $\{l_{ij}, l_{jk}, l_{ik}\}$ satisfy the triangle inequalities for any face $\triangle ijk$,
$r$ is called an inversive distance circle packing metric. Denote

\begin{equation*}
\Omega_{ijk}=\{(r_i, r_j, r_k)|l_{ij}<l_{ik}+l_{kj}, l_{ik}<l_{ij}+l_{jk}, l_{jk}<l_{ji}+l_{ik} \}
\end{equation*}
and
\begin{equation*}
\Omega=\cap_{\Delta ijk\in F}\Omega_{ijk}=\{r\in \mathds{R}_{>0}^N|l_{ij}<l_{ik}+l_{kj}, l_{ik}<l_{ij}+l_{jk}, l_{jk}<l_{ji}+l_{ik}, \forall \Delta ijk\in F\}.
\end{equation*}
If $r\in \Omega$, the triangulated surface could be obtained by gluing many triangles along their edges coherently and there will
be singularities at the vertices.  The classical discrete Gauss curvature $K$ is introduced to describe the singularity, which is
defined as the angle deficit at a vertex, i.e. $K_i=2\pi-\sum_{\triangle ijk\in F}\theta_i^{jk}$,
where $\theta_i^{jk}$ is the inner angle of $\triangle ijk$ at $i$. In this paper, we introduce
a new combinatorial curvature, which is defined as
$$R_i=\frac{K_i}{r_i^2}$$
for Euclidean background geometry and
$$R_i=\frac{K_i}{\tanh^2 \frac{r_i}{2}}$$
for hyperbolic background geometry.
For this combinatorial curvature, we prove the following main theorem on global rigidity.
\begin{theorem}\label{rigidity for R}
Given a weighted closed triangulated surface $(M,\mathcal{T}, I)$ with $I\geq0$ and $\chi(M)\leq0$.
$\overline{R}\in C(V)$ is a given function on $V$.
\begin{description}
  \item[(1)] If $\overline{R}\equiv0$, then there exists at most one Euclidean inversive distance circle packing metric $\overline{r}\in \Omega$
    up to scaling such that its $R$-curvature is 0. If $\overline{R}\leq0$ and $\overline{R}\not\equiv 0$, then
      there exists at most one Euclidean inversive distance circle packing metric $\overline{r}\in \Omega$
        such that its $R$-curvature is the given function $\overline{R}$.
  \item[(2)] If $\overline{R}\leq 0$, then there exists at most one hyperbolic inversive distance circle packing metric
    $\overline{r}\in \Omega$ such that its $R$-curvature is the given function $\overline{R}$.
\end{description}
\end{theorem}

We further give a counterexample, i.e. Example \ref{counter example},
which shows that there is no rigidity for triangulated surfaces with positive Euler number.
So the result in Theorem \ref{rigidity for R} is sharp.

Following \cite{GX4}, we generalize the combinatorial curvature $R_i$ to the $\alpha$-curvature $R_{\alpha}$,
which is defined as
$R_{\alpha, i}=\frac{K_i}{r_i^\alpha}$ for Euclidean background
geometry and $R_{\alpha, i}=\frac{K_i}{\tanh^\alpha\frac{r_i}{2}}$ for hyperbolic background geometry.
For the $\alpha$-curvature $R_{\alpha}$, we have the following generalized global rigidity.
\begin{theorem}\label{main rigidity}
Given a weighted closed triangulated surface $(M,\mathcal{T}, I)$ with $I\geq0$ and $\alpha\chi(M)\leq0$.
$\overline{R}\in C(V)$ is a given function on $V$.
\begin{description}
  \item[(1)]
  If $\alpha \overline{R}\equiv0$, then there exists at most one Euclidean inversive distance
  circle packing metric $\overline{r}\in \Omega$ up to scaling
  with  $\alpha$-curvature $\overline{R}$;
  If $\alpha\overline{R}\leq 0$ and $\alpha\overline{R}\not\equiv 0$, then there
  exists at most one Euclidean inversive distance circle packing metric $\overline{r}\in \Omega$
  with  $\alpha$-curvature $\overline{R}$.
  \item[(2)]
  If $\alpha\overline{R}\leq 0$, then there
  exists at most one hyperbolic inversive distance circle packing metric $\overline{r}\in \Omega$
  with  $\alpha$-curvature $\overline{R}$.
\end{description}
\end{theorem}

The result in Theorem \ref{main rigidity} is also sharp.
In the special case of $\alpha=0$, this is Bowers and Stephenson's conjecture on the global rigidity
of the classical combinatorial Gauss curvature $K$ for inversive distance circle packing metrics, which was proved by
Guo \cite{Guo} and Luo \cite{L2}. So Theorem \ref{main rigidity} in fact partially solves a generalized Bowers-Stephenson conjecture.

To study the constant curvature problem of the curvature $R$ for the Euclidean background geometry,
we introduce a combinatorial Ricci flow, defined as
\begin{equation}\label{Euclidean Ricci flow introduction}
\frac{dg_i}{dt}=(R_{av}-R_i)g_i,
\end{equation}
where $g_i=r_i^2$ and $R_{av}=\frac{2\pi\chi(M)}{||r||_2^2}$.
We have the following main result on the the combinatorial Euclidean Ricci flow (\ref{Euclidean Ricci flow introduction}) and
existence of constant $R$-curvature metrics.

\begin{theorem}\label{constant curvature problem introduction}
Given a weighted closed triangulated surface $(M,\mathcal{T}, I)$ with $I\geq0$.
\begin{description}
  \item[(1)] Along the Euclidean Ricci flow (\ref{Euclidean Ricci flow introduction}),
the curvature $R$ evolves according to
$$\frac{dR_{i}}{dt}=\Delta R_{i}+R_{i}(R_{i}-R_{av}),$$
where the Laplace operator is defined to be
\begin{equation*}
\Delta f_i=\frac{1}{r_i^2}\sum_{j=1}^N \left(-\frac{\partial K_i}{\partial u_j}\right) f_j
=\frac{1}{r_i^2}\sum_{j\sim i} \left(-\frac{\partial K_i}{\partial u_j}\right) (f_j-f_i)
\end{equation*}
with $u_i=\ln r_i^2$ for $f\in C(V)$.
  \item[(2)] If the solution of the Euclidean combinatorial Ricci flow (\ref{Euclidean Ricci flow introduction}) stays in $\Omega$ and
converges to an inversive distance circle packing metric $r^*\in \Omega$, then there exists a constant $R$-curvature
metric in $\Omega$ and $r^*$ is such one.
  \item[(3)] If $\chi(M)\leq 0$ and there exists a Euclidean inversive distance circle packing metric $r^*\in \Omega$ with constant $R$-curvature,
then the solution of the Euclidean combinatorial Ricci flow (\ref{Euclidean Ricci flow introduction}) develops no essential singularities
in finite time and at time infinity; Furthermore, if the solution of (\ref{Euclidean Ricci flow introduction}) develops no
removable singularities in finite time, then the solution of (\ref{Euclidean Ricci flow introduction}) exists for all time,
converges exponentially fast to the constant $R$-curvature metric $r^*$ and does not develop removable
singularities at time infinity.
  \item[(4)] If $\chi(M)\leq 0$ and there exists a Euclidean inversive distance
circle packing metric $r^*\in \Omega$ with constant $R$-curvature, then the solution of the extended combinatorial
Ricci flow
\begin{equation}\label{extended Ricci flow introduction}
\frac{dr_i}{dt}=(R_{av}-\widetilde{R}_{i})g_i
\end{equation}
exists for all time and converges exponentially fast to $r^*$.
\end{description}
\end{theorem}

The definition of essential singularity removable singularity is given in Definition \ref{definition of singularity}
and $\widetilde{R}$ is an extension of the $R$-curvature. The existence of constant $R$-curvature metric is in fact
equivalent to the convergence of the extended Ricci flow (\ref{extended Ricci flow introduction}).

For the hyperbolic background geometry, the most important constant curvature problem is the zero curvature problem, which
frequently appears in engineering. Following \cite{GX3}, we treat it as a special case of the prescribing curvature problem.
Interestingly, we found that the results of the prescribed curvature problems
for Euclidean and hyperbolic background geometry are similar,
so we state them together in a unified form here.

\begin{theorem}\label{prescribing curvature problem introduction}
Given a closed triangulated surface $(M, \mathcal{T})$ with inversive distance $I\geq 0$.
Given a function $\overline{R}\in C(V)$, the modified combinatorial Ricci flow is defined to be
\begin{equation}\label{modifed Ricci flow introduction}
\frac{dg_i}{dt}=(\overline{R}_i-R_{i})g_i,
\end{equation}
where $g_i=r_i^2$ for the Euclidean background geometry and $g_i=\tanh^2\frac{r_i}{2}$ for the hyperbolic background geometry.
\begin{description}
  \item[(1)] If the solution of the modified Ricci flow (\ref{modifed Ricci flow introduction})
stays in $\Omega$ and converges to $\overline{r}\in \Omega$, then we have $\overline{R}$ is admissible and $R(\overline{r})=\overline{R}$.
  \item[(2)] Suppose $\overline{R}\leq 0$ and $\overline{R}$ is admissible with $R(\overline{r})=\overline{R}$,
then the modified Ricci flow (\ref{modifed Ricci flow introduction})
develops no essential singularity in finite time and at time infinity.
Furthermore, if the solution of the Ricci flow (\ref{modifed Ricci flow introduction})
develops no removable singularities at finite time, then the solution of (\ref{modifed Ricci flow introduction})
exists for all time, converges exponentially fast to $\overline{r}$ and does not develop removable singularities at time infinity.
  \item[(3)] If
  $\overline{R}\in C(V)$ is admissible,
  $R(\overline{r})=\overline{R}$ with $\overline{r}\in \Omega$ and $\overline{R}\leq 0$,
  then any solution of the extended combinatorial Ricci flow
  $$\frac{dg_i}{dt}=(\overline{R}_i-\widetilde{R}_{i})g_i$$
  exists for all time and converges exponentially fast to $\overline{r}$.
\end{description}
\end{theorem}

We can also introduce a combinatorial $\alpha$-Ricci flow to study the constant and prescribed curvature problem for the
$\alpha$-curvature $R_{\alpha}$. The results are parallel to Theorem \ref{constant curvature problem introduction} and Theorem
\ref{prescribing curvature problem introduction} respectively.
The precise statements of the results are given in Theorem \ref{constant alpha curvature problem} and
Theorem \ref{prescribling alpha curvature problem} respectively, so we do not list them here.

The paper is organized as follows.
In Section \ref{section defnition of combinatorial Gauss curvature},
we establish the framework and recall some basic facts on inversive distance circle packing metrics and
the classical Gauss curvature,
then we introduce the new definition of combinatorial curvature for triangulated surfaces
with inversive distance circle packing metrics.

In Section \ref{section rigidity of the combinatorial curvature}, we give a detailed proof of Theorem \ref{rigidity for R}.
We first prove the local rigidity of the new combinatorial curvature $R$ with respect to the inversive distance circle
packing metrics by variational principles, i.e. Theorem \ref{local rigidity}.
Then we prove the global rigidity of the new combinatorial curvature $r$
with respect to the inversive distance circle packing metrics, i.e. Theorem \ref{theorem global rigidity},
using the extension of convex functions by constants in \cite{BPS, L2}.
Theorem \ref{local rigidity} and Theorem \ref{theorem global rigidity} together generate Theorem \ref{rigidity for R}.
We further give a counterexample,
i.e. Example \ref{counter example}, which shows that the global rigidity in Theorem \ref{theorem global rigidity} is sharp.

In Section \ref{combinatorial Ricci flow}, we first pose the combinatorial Yamabe problem of the combinatorial curvature $R$
with respect to the inversive distance circle packing metrics and then we introduce the combinatorial Ricci flow
to study the combinatorial Yamabe problem.
In subsection \ref{subsection Euclidean combinatorial Ricci flow},
we study the behavior of the Euclidean combinatorial Ricci flow. We first derive the evolution of $R$
along the Euclidean combinatorial Ricci flow (\ref{Euclidean Ricci flow introduction}) in Lemma \ref{evolution of R along Euclidean Ricci flow},
then we study the necessary condition for the convergence of the flow (\ref{Euclidean Ricci flow introduction}) and obtain
Proposition \ref{necessary for existence of cont curv metric}.
To study the long time behavior of the flow , we introduce the notion of singularities in Definition \ref{definition of singularity}.
Then we get Theorem \ref{nonexistence of essential singularity}, Theorem \ref{convergence under nonexistence of removable singularity}
and Theorem \ref{convergence of extended flow} with the aid of Ricci potential function and extended combinatorial Ricci flow.
Lemma \ref{evolution of R along Euclidean Ricci flow}, Proposition \ref{necessary for existence of cont curv metric},
Theorem \ref{nonexistence of essential singularity}, \ref{convergence under nonexistence of removable singularity}
and \ref{convergence of extended flow} together generate Theorem \ref{constant curvature problem introduction}.
We also use the Euclidean combinatorial Ricci flow to study the prescribed curvature problem and
obtain Theorem \ref{Euclidean prescribing curvature problem}. In subsection \ref{subsection hyperbolic combinatorial Ricci flow},
we study the hyperbolic prescribed curvature problem using the hyperbolic combinatorial Ricci flow.
We obtain the necessary condition for the convergence of the hyperbolic flow in
Proposition \ref{necessary for the convergence of the midified hyper flow} and
the main result on hyperbolic prescribed curvature problem in Theorem \ref{hyperbolic prescribing curvature problem}.
Theorem \ref{Euclidean prescribing curvature problem}, Proposition \ref{necessary for the convergence of the midified hyper flow}
and Theorem \ref{hyperbolic prescribing curvature problem} together generate Theorem \ref{prescribing curvature problem introduction}.

In Section \ref{section alpha curvature and alpha curvature flows}, we generalize the definition of $R$-curvature to $\alpha$-curvature
and then study the rigidity of the curvature, the constant curvature problem and the prescribed curvature problem.
We obtain Theorem \ref{main rigidity}, Theorem \ref{constant alpha curvature problem} and
Theorem \ref{prescribling alpha curvature problem} in this section.

\section{Definition of combinatorial Gauss curvature}\label{section defnition of combinatorial Gauss curvature}
In this section, we give the preliminaries that are needed in the paper, including the notations,
the definition of inversive distance circle packing metrics. Then we introduce the new combinatorial Gauss
curvature for inversive distance circle packing metrics.

Suppose $M$ is a closed surface with a triangulation $\mathcal{T}=\{V,E,F\}$,
where $V,E,F$ represent the sets of vertices, edges and faces respectively.
Let $I: E\rightarrow [0,+\infty)$ be a function assigning each edge $\{ij\}$ a weight $I_{ij}\in [0,+\infty)$,
which is denoted as $I\geq 0$ in the paper.
The triple $(M, \mathcal{T}, I)$ will be referred to as a weighted triangulation of $M$ in the following.
All the vertices are ordered one by one, marked by $v_1, \cdots, v_N$, where $N$
is the number of vertices, and we often use $i$ to denote the vertex $v_i$ for simplicity in the following.
We use $i\sim j$ to denote that the
vertices $i$ and $j$ are adjacent if there is an edge $\{ij\}\in E$ with $i$, $j$ as end points.
Throughout this paper, all functions $f: V\rightarrow \mathds{R}$ will be regarded as column
vectors in $\mathds{R}^N$ and $f_i$ is the value of $f$ at $i$. And we use $C(V)$ to denote the set of functions
defined on $V$.

Each map $r:V\rightarrow (0,+\infty)$ is a circle packing, which could be taken as the radius $r_i$ of a circle
attached to the vertex $i$.
Given $(M, \mathcal{T}, I)$, we attach each edge $\{ij\}$ the length
\begin{equation}\label{definition of length of edge for Euclidean background}
l_{ij}=\sqrt{r_i^2+r_j^2+2r_ir_jI_{ij}}
\end{equation}
for Euclidean background geometry and
\begin{equation}\label{definition of length of edge for hyperbolic background}
l_{ij}=\cosh^{-1}(\cosh (r_i)\cosh(r_j)+I_{ij}\sinh(r_i)\sinh(r_j))
\end{equation}
for hyperbolic background geometry.
Then $I_{ij}$ is the inversive distance of the two circles centered at $v_i$ and $v_j$ with radii $r_i$ and $r_j$ respectively.
If $I_{ij}\in[0,1]$, the two circles intersect at an acute angle and we can take $I_{ij}=\cos \Phi_{ij}$
with $\Phi_{ij}\in [0, \frac{\pi}{2}]$ and then the inversive distance circle packing is reduced to
Thurston's circle packing.
However, for general weight $I_{ij}\in [0,+\infty)$,  in order that the lengths $l_{ij}, l_{jk}, l_{ik}$ for a
face $\Delta ijk\in F$ satisfy the triangle inequalities, the admissible space of the radius
\begin{equation}\label{admissible space}
\Omega_{ijk}=\{(r_i, r_j, r_k)|l_{ij}<l_{ik}+l_{kj}, l_{ik}<l_{ij}+l_{jk}, l_{jk}<l_{ji}+l_{ik} \}
\end{equation}
is not $\mathds{R}_{>0}^3$.
In fact, it is proved \cite{Guo} that the admissible space $\Omega_{ijk}$ is a simply connected open subset of $\mathds{R}_{>0}^3$.
Note that the set $\Omega_{ijk}$ maybe not convex.
Set
\begin{equation}
\Omega=\cap_{\Delta ijk\in F}\Omega_{ijk}=\{r\in \mathds{R}_{>0}^N|l_{ij}<l_{ik}+l_{kj}, l_{ik}<l_{ij}+l_{jk}, l_{jk}<l_{ji}+l_{ik}, \forall \Delta ijk\in F\}
\end{equation}
to be the space of admissible radius function. $\Omega$ is obviously an open subset of $\mathds{R}^N_{>0}$.
Every $r\in \Omega$ is called a circle packing metric.
For more information on inversive distance circle packing metrics, the readers can refer to  Stephenson \cite{St},
Bowers and Hurdal \cite{BH} and Guo \cite{Guo}.

Now suppose that for each face $\Delta ijk\in F$, the triangle inequalities are satisfied, i.e. $r\in \Omega$,
then the weighted triangulated surface
$(M, \mathcal{T}, I)$ could be taken as gluing many triangles along the edges coherently, which produces
a cone metric on the triangulated surface with singularities at the vertices.
To describe the singularity at the vertex $i$, the classical discrete Gauss curvature is introduced, which is defined as
\begin{equation}\label{classical Gauss curv}
K_i=2\pi-\sum_{\triangle ijk \in F}\theta_i^{jk},
\end{equation}
where the sum is taken over all the triangles with $i$ as one of its vertices. Discrete Gaussian curvature $K_i$ satisfies the following discrete version of Gauss-Bonnet formula \cite{CL1}
\begin{equation}\label{Gauss-Bonnet without weight}
\sum_{i\in V}K_i=2\pi \chi(M)-\lambda Area(M),
\end{equation}
where $\lambda=0, -1$ in the case of Euclidean and hyperbolic background geometry respectively.
Guo \cite{Guo} further proved the following property for the classical combinatorial Gauss curvature.
\begin{lemma}\label{positivity of L}
(Guo \cite{Guo}) Given $(M, \mathcal{T}, I)$ with inversive distance $I\geq0$. Set $u_i=\ln r_i^2$ for the Euclidean background geometry,
$u_i=\ln \tanh^2 \frac{r_i}{2}$ for the hyperbolic background geometry and $L=\frac{\partial (K_1,\cdots, K_N)}{\partial(u_1,\cdots,u_N)}$. Then
\begin{description}
  \item[(1)] $L$ is symmetric and positive semi-definite with rank $N-1$ and kernel $\{t\textbf{1}|t\in \mathds{R}\}$ on $\Omega$ for the Euclidean background geometry;
  \item[(2)] $L$ is symmetric and positive definite on $\Omega$ for the hyperbolic background geometry.
\end{description}
\end{lemma}
Lemma \ref{positivity of L} implies the local rigidity of the inversive distance circle packing metric obtained in \cite{Guo}.
Based on Guo's work on local rigidity \cite{Guo} and Bobenko, Pinkall and Springborn's work \cite{BPS} on convex extension of functions,
Luo \cite{L2} proved the following global rigidity for inversive distance circle packing metrics.

\begin{theorem}
(Luo \cite{L2}) Given a weighted closed triangulated surface $(M, \mathcal{T}, I)$ with inversive distance $I\geq 0$.
\begin{description}
  \item[(1)] A Euclidean inversive distance circle packing metric on $(M, \mathcal{T}, I)$ is uniquely
  determined by its combinatorial Gauss curvature $K$ up to scaling.
  \item[(2)] A hyperbolic inversive distance circle packing metric on $(M, \mathcal{T}, I)$ is uniquely
  determined by its combinatorial Gauss curvature $K$.
\end{description}
\end{theorem}

As noted in \cite{GX4}, the classical definition of combinatorial Gauss curvature $K_i$ with Euclidean background geometry
in (\ref{classical Gauss curv}) has two disadvantages. The first is that the classical combinatorial Gauss curvature is scaling
invariant, i.e. $K_i(\lambda r)=K_i(r)$ for any $\lambda>0$; The second is that, as the triangulated surfaces approximate a smooth surface,
the classical combinatorial curvature $K_i$ could not approximate the smooth Gauss curvature, as we obviously have $K_i$ tends zero.
Motivated by the two disadvantages, we introduce a new combinatorial Gauss curvature for triangulated surfaces with Thurston's circle
packing metrics in \cite{GX4} and we can generalize the curvature to the case of inversive distance circle packing metrics here.

For the following applications, we always set
\begin{equation}\label{s_i}
s_i(r)=\left\{
         \begin{array}{ll}
           r_i, & \hbox{Euclidean background geometry} \\
           \tanh\frac{r_i}{2}, & \hbox{hyperbolic background geometry}
         \end{array}
       \right.
\end{equation}
in this paper.
Then we can generalize the definition of combinatorial Gauss curvature to triangulated surfaces
with inversive distance circle packing metrics as follows.

\begin{definition}\label{definition of comb curv with inversive dist}
Given a triangulated surface $(M, \mathcal{T})$ with inversive distance $I\geq 0$
and a circle packing metric $r\in \Omega$, the combinatorial Gauss curvature at the vertex $i$ is defined
to be
\begin{equation}\label{definition of Gauss curv}
R_i=\frac{K_i}{s_i^2},
\end{equation}
where $K_i$ is the classical combinatorial Gauss curvature at $i$ given by (\ref{classical Gauss curv}) and $s_i$ is given by (\ref{s_i}).
\end{definition}

As the inversive distance generalizes Thurston's intersection angle, the Definition \ref{definition of comb curv with inversive dist}
of combinatorial Gauss curvature naturally generalizes the definition of combinatorial Gauss curvature in \cite{GX4,GX3}.
It is obvious that the curvature $R_i$ is an elementary function of $r$ and then obviously a smooth function of the radius function $r$.
In the following, we often refer to the combinatorial Gauss curvature
in Definition \ref{definition of comb curv with inversive dist} as $R$-curvature for short.

As the Riemannian metric tensor is a positive definite symmetric 2-tensor in Riemannian geometry, we can take $g_i=r_i^2$ as the
analogue of the Riemannian metric in the combinatorial setting with Euclidean background geometry. Then we have
$$R_i(\lambda g)=\lambda^{-1}R_i(g)$$
for any constant $\lambda>0$. Furthermore, Example 1 in \cite{GX4} shows that the combinatorial Gauss curvature $R_i$ could
approximate the smooth Gauss curvature up to a uniform constant $\pi$ for the Euclidean background geometry.
Both facts indicate that $R_i$
is a good candidate as the combinatorial Gauss curvature on triangulated surfaces with inversive distance circle packing metrics.

If we take the measure at the vertex $v_i$ as $\mu_i=r_i^2$ for the Euclidean background geometry,
then we have
$\int_{M}fd\mu=\sum_{i=1}^N f_ir_i^2$
for any $f\in C(V)$. Using this discrete measure, we have the following combinatorial Gauss-Bonnet formula as in \cite{GX4}
\begin{equation}\label{Gauss-Bonnet with weight}
\int_MRd\mu=\sum_iR_id\mu_i=2\pi\chi(M).
\end{equation}
In this sense, the average curvature is
\begin{equation}\label{avarage curvature}
R_{av}=\frac{\int_MRd\mu}{\int_Md\mu}=\frac{2\pi\chi(M)}{||r||_2^2},
\end{equation}
where $||r||_2^2=\sum_{i=1}^Nr_i^2$ is the total measure of $M$ with respect to $\mu$.

For the hyperbolic background geometry, if we take $\mu_i=\tanh^2\frac{r_i}{2}$, we also have the following
Gauss-Bonnet formula
\begin{equation}\label{hyperbolic Gauss-Bonnet formula}
\int_MRd\mu=\sum_iR_id\mu_i=2\pi\chi(M)+Area(M).
\end{equation}

\section{Rigidity of the combinatorial curvature}\label{section rigidity of the combinatorial curvature}

In this section, we prove the global rigidity of the $R$-curvature in the sense
that the $R$-curvature is uniquely determined by the inversive distance
circle packing metrics.
We first prove the local rigidity of the $R$-curvature with respect to the inversive distance circle packing metric $r$.
Before that, we state the following results obtained by Guo \cite{Guo} and Luo \cite{L2}.

\begin{lemma}\label{extension lemma}
(\textbf{Guo-Luo})
Given a weighted triangulated surface $(M, \mathcal{T})$ with inversive distance $I\geq 0$,
suppose $\Delta ijk\in F$ and $\theta_i$ is the inner angle of the triangle at the vertex $i$,
then we have
\begin{description}
  \item[(1)] the function
\begin{equation}\label{Guo's energy function}
F_{ijk}(u)=\int_{u_0}^u \theta_idu_i+\theta_jdu_j+\theta_kdu_k
\end{equation}
is well defined on
$$\mathcal{U}_{ijk}=\{(u_i,u_j,u_k)|u_i=\ln s_i^2, l_{ij}<l_{ik}+l_{kj}, l_{ik}<l_{ij}+l_{jk}, l_{jk}<l_{ji}+l_{ik}\},$$
where $u_0\in \mathcal{U}_{ijk}$;
  \item[(2)] $F_{ijk}(u)$ is strictly concave on $\mathcal{U}_{ijk}\cap \{u_i+u_j+u_k=0\}$ for the Euclidean background geometry and strictly concave on $\mathcal{U}_{ijk}$
for the hyperbolic background geometry;
  \item[(3)] the function $F_{ijk}(u)$ could be extended by constant to a $C^1$ smooth concave function
\begin{equation}\label{extension of F_{ijk}}
\widetilde{F}_{ijk}(u)=\int_{u_0}^u \widetilde{\theta}_idu_i+\widetilde{\theta}_jdu_j+\widetilde{\theta}_kdu_k
\end{equation}
defined on $u\in \mathds{R}^3$ for the Euclidean background geometry and $u\in \mathds{R}^N_{<0}$ for the hyperbolic background geometry,
where the extension $\widetilde{\theta}_i$ of $\theta_i$ by constant is defined to be
$\widetilde{\theta}_i=\pi$ when $l_{jk}\geq l_{ik}+l_{ij}$ and $\widetilde{\theta}_i=0$
when $l_{ik}\geq l_{jk}+l_{ij}$ or $l_{ij}\geq l_{jk}+l_{ik}$.
\end{description}

\end{lemma}

Using this result, we have the following local rigidity for the $R$-curvature.

\begin{proposition}\label{local rigidity}
Given a triangulated surface $(M, \mathcal{T})$ with inversive distance $I\geq 0$ and $\chi(M)\leq 0$.
$\overline{R}\in C(V)$ is a given function defined on $(M, \mathcal{T})$.
\begin{description}
 \item[(1)]  In the case of Euclidean background geometry,
  if $\overline{R}\equiv0$, then there locally exists at most one inversive distance circle packing metric $\overline{r}\in \Omega$
with $R$-curvature $\overline{R}$ up to scaling;
if $\overline{R}\leq 0$ and $\overline{R}\not\equiv 0$, then
there locally exists at most one inversive distance circle packing metric $\overline{r}\in \Omega$
with $R$-curvature $\overline{R}$.
  \item[(2)] In the case of hyperbolic background geometry, if $\overline{R}\leq 0$,
then there locally exists at most one inversive distance circle packing metric $\overline{r}\in \Omega$
with combinatorial $R$-curvature $\overline{R}$.
\end{description}
\end{proposition}

\proof
We first prove the case of Euclidean background geometry.
Suppose $r_0\in \Omega$ is an inversive distance circle packing metric and $u_{0,i}=\ln r_{0,i}^2$.
Define a Ricci potential function
\begin{equation}\label{Euclidean energy function}
\begin{aligned}
F(u)=-\sum_{\Delta ijk\in F}F_{ijk}(u)+\int_{u_0}^u\sum_{i=1}^N(2\pi-\overline{R}_ir_i^2)du_i,
\end{aligned}
\end{equation}
where $F_{ijk}(u)$ is given by (\ref{Guo's energy function}) with $s_i=r_i$.
Note that the second term $\int_{u_0}^u\sum_{i=1}^N(2\pi-\overline{R}_ir_i^2)du_i$
in (\ref{Euclidean energy function}) is independent of the paths from $u_0$ to $u$.
Combined with Lemma \ref{extension lemma},
this implies the Ricci potential functional $F(u)$ is well-defined on $\mathcal{U}=\cap_{\triangle ijk\in F}\mathcal{U}_{ijk}$ and
$$\nabla_{u_i}F=-\sum_{\Delta ijk\in F}\theta_i^{jk}+2\pi-\overline{R}_ir_i^2=K_i-\overline{R}_ir_i^2.$$
Then $r$ is an inversive distance circle packing metric with combinatorial curvature $\overline{R}$ if and only if
the corresponding $u$ is a critical point of $F(u)$.
Note that
\begin{equation}
\operatorname{Hess}_u F= L-\left(
              \begin{array}{ccc}
                \overline{R}_1r_i^2 &   &   \\
                  & \ddots &   \\
                  &   & \overline{R}_Nr_N^2 \\
              \end{array}
            \right),
\end{equation}
we have $\operatorname{Hess}_u F$ is symmetric and positive semi-definite by Lemma \ref{positivity of L}, if $\overline{R}\leq 0$.
In the case of $\overline{R}\equiv0$, the kernel of $\operatorname{Hess}_u F$ is $\{t\textbf{1}|t\in \mathds{R}\}$.
While in the case of $\overline{R}\leq 0$ and $\overline{R}\not\equiv 0$, we have
$\operatorname{Hess}_u F$ is strictly positive definite. For a strictly convex function, its gradient is locally injective.
This implies the result for the Euclidean case, i.e. the first part of the proposition.

The case for hyperbolic background geometry could be proved similarly using the following Ricci potential function
\begin{equation*}
\begin{aligned}
F(u)=-\sum_{\Delta ijk\in F}F_{ijk}(u)+\int_{u_0}^u\sum_{i=1}^N(2\pi-\overline{R}_i\tanh^2\frac{r_i}{2})du_i.
\end{aligned}
\end{equation*}
As the proof is almost the same as that for the Euclidean background geometry, we omit the details here.
\qed

\begin{remark}
If $\Omega$ is simply connected, the Ricci potential function $F(u)$ could be written in a simple form
\begin{equation}\label{good form of ricci potential F}
F(u)=\int_{u_0}^u\sum_{i=1}^N(K_i-\overline{R}_is_i^2)du_i.
\end{equation}
Especially, for the case of $I\in [0, 1]$, the inversive distance circle packing metric is just Thurston's circle packing
metric, $\Omega=\mathds{R}^N_{>0}$ and the Ricci potential function $F(u)$
could be written in the form of (\ref{good form of ricci potential F}) \cite{GX4, GX3}.
\end{remark}

In the following of this section, we focus on the global rigidity of the combinatorial curvature $R$ with respect to
the inversive distance circle packing metrics.
Using the extension method, we can prove the following main theorem on global rigidity.

\begin{theorem}\label{theorem global rigidity}
Given a triangulated surface $(M, \mathcal{T})$ with inversive distance $I\geq 0$ and $\chi(M)\leq 0$.
$\overline{R}\in C(V)$ is a given function defined on $(M, \mathcal{T})$.
\begin{description}
  \item[(1)]  In the case of Euclidean background geometry,
  if $\overline{R}\equiv0$, there exists at most one inversive distance circle packing metric $\overline{r}\in \Omega$
  with combinatorial $R$-curvature $\overline{R}$ up to scaling;
  if $\overline{R}\leq 0$ and $\overline{R}\not\equiv0$,
  there exists at most one inversive distance circle packing metric $\overline{r}\in \Omega$ with combinatorial $R$-curvature $\overline{R}$.
  \item[(2)] In the case of hyperbolic background geometry, if $\overline{R}\leq 0$,
  there exists at most one inversive distance circle packing metric $\overline{r}\in \Omega$
  with combinatorial $R$-curvature $\overline{R}$.
\end{description}
\end{theorem}

\proof
We only prove the case of Euclidean background geometry and the case of hyperbolic background geometry could be proved similarly.

The Ricci potential function $F(u)$ could be extended from $\mathcal{U}$ to the whole space $\mathds{R}^{N}$, where
$\mathcal{U}$ is the image of $\Omega$ under the map $u_i=\ln r_i^2$.
Suppose $r_0\in \Omega$ is an inversive distance circle packing metric and $u_0=\ln r_{0, i}^2$.
Then the Ricci energy  $F(u)$ defined by $\overline{R}$ is
\begin{equation}\label{F extend}
\begin{aligned}
F(u)=-\sum_{\Delta ijk\in F}F_{ijk}(u)+\int_{u_0}^u\sum_{i=1}^N(2\pi-\overline{R}_ir_i^2)du_i.
\end{aligned}
\end{equation}
Note that the function $F_{ijk}(u)$ defined on $\mathcal{U}_{ijk}$
could be extended to $\widetilde{F}_{ijk}(u)$ defined by (\ref{extension of F_{ijk}}) on $\mathds{R}^N$ by Lemma \ref{extension lemma}
and the second term $\int_{u_0}^u\sum_{i=1}^N(2\pi-\overline{R}_ir_i^2)du_i$ in the second line of (\ref{F extend})
can be naturally defined on $\mathds{R}^N$, then we have the following extension $\widetilde{F}(u)$ defined on $\mathds{R}^N$
of the Ricci potential function $F(u)$
$$\widetilde{F}(u)=-\sum_{\triangle ijk\in F}\widetilde{F}_{ijk}(u)+\int_{u_0}^u\sum_{i=1}^N(2\pi-\overline{R}_ir_i^2)du_i.$$
As $\widetilde{F}_{ijk}(u)$ is $C^1$-smooth concave by Lemma \ref{extension lemma} and $\int_{u_0}^u\sum_{i=1}^N(2\pi-\overline{R}_ir_i^2)du_i$
is a well-defined convex function on $\mathds{R}^N$ for $\overline{R}\leq 0$,
we have $\widetilde{F}(u)$ is a $C^1$-smooth convex function on $\mathds{R}^N$. Furthermore,
\begin{equation*}
\nabla_{u_i}\widetilde{F}=-\sum_{\triangle ijk\in F}\widetilde{\theta}_i+2\pi-\overline{R}_ir_i^2=\widetilde{K}_i-\overline{R}_ir_i^2,
\end{equation*}
where $\widetilde{K}_i=2\pi-\sum_{\triangle ijk\in F}\widetilde{\theta}_i$.
By the proof of Proposition \ref{local rigidity}, we have $\widetilde{F}(u)$ is convex on $\mathds{R}^N$ and strictly convex on $\mathcal{U}\cap \{\sum_{i=1}^Nu_i=0\}$ for $\overline{R}\equiv0$.
Similarly, we have $\widetilde{F}(u)$ is convex on $\mathds{R}^N$ and strictly convex on
$\mathcal{U}$ for $\overline{R}\leq0$ and $\overline{R}\not\equiv0$.

If there are two different inversive distance circle packing metrics $\overline{r}_{A}, \overline{r}_{B}\in \Omega$ with the same combinatorial
$R$-curvature $\overline{R}$, then
$\overline{u}_A=\ln \overline{r}_{A}^2\in \mathcal{U}$, $\overline{u}_B=\ln \overline{r}_{B}^2\in \mathcal{U}$
are both critical points of the extended Ricci potential $\widetilde{F}(u)$.
It follows that
$$\nabla \widetilde{F}(\overline{u}_A)=\nabla \widetilde{F}(\overline{u}_B)=0.$$
Set
\begin{equation*}
\begin{aligned}
f(t)=&\widetilde{F}((1-t)\overline{u}_A+t\overline{u}_B)\\
=&\sum_{\triangle ijk\in F}f_{ijk}(t)+\int_{u_0}^{(1-t)\overline{u}_A+t\overline{u}_B}\sum_{i=1}^N(2\pi-\overline{R}_ir_i^2)du_i,
\end{aligned}
\end{equation*}
where
$$f_{ijk}(t)=-\widetilde{F}_{ijk}((1-t)\overline{u}_A+t\overline{u}_B).$$
Then $f(t)$ is a $C^1$ convex function on $[0, 1]$ and $f'(0)=f'(1)=0$, which implies that $f'(t)\equiv 0$ on $[0, 1]$.
Note that $\overline{u}_A$ belongs to the open set $\mathcal{U}$,
there exists $\epsilon>0$ such that $(1-t)\overline{u}_A+t\overline{u}_B\in \mathcal{U}$ for $t\in [0, \epsilon]$.
So $f(t)$ is smooth on $[0, \epsilon]$.

In the case of $\overline{R}\leq 0$ and $\overline{R}\not\equiv0$, the strict convexity of $\widetilde{F}(u)$
on $\mathcal{U}$ implies that $f(t)$ is strictly convex on $[0, \epsilon]$ and $f'(t)$ is a strictly increasing function on $[0, \epsilon]$.
Then $f'(0)=0$ implies $f'(\epsilon)>0$, which contradicts $f'(t)\equiv 0$ on $[0, 1]$. So there exists at most one inversive circle packing
metric with combinatorial $R$-curvature $\overline{R}$.

For the case of $\overline{R}\equiv0$, we have $f(t)$ is $C^1$ convex on $[0, 1]$ and smooth on $[0, \epsilon]$.
$f'(t)\equiv 0$ on $[0, 1]$ implies that $f''(t)\equiv 0$ on $[0, \epsilon]$.
Note that, for $t\in [0, \epsilon]$,
\begin{equation*}
\begin{aligned}
f''(t)=(\overline{u}_A-\overline{u}_B) L  (\overline{u}_A-\overline{u}_B)^T.
\end{aligned}
\end{equation*}
By Lemma \ref{positivity of L}, we have $\overline{u}_A-\overline{u}_B=c\textbf{1}$ for some constant $c\in \mathds{R}$, which
implies that $\overline{r}_A=e^{c/2}\overline{r}_B$. So there exists at most one inversive distance circle packing metric
with zero $R$-curvature up to scaling.
\qed

As a corollary of Theorem \ref{theorem global rigidity}, we have the following global rigidity for the inversive distance circle packing metrics with constant combinatorial $R$-curvature.
\begin{corollary}
Given a closed triangulated surface $(M, \mathcal{T})$ with inversive distance $I\geq 0$ and $\chi(M)\leq 0$.
\begin{description}
  \item[(1)]  In the case of Euclidean background geometry,
  for $\chi(M)=0$, there exists at most one inversive distance circle packing metric
  with constant combinatorial $R$-curvature $0$ up to scaling.
  For $\chi(M)<0$,
  there exists at most one inversive distance circle packing metric with negative constant combinatorial $R$-curvature $c\in \mathds{R}$.
  \item[(2)] In the case of hyperbolic background geometry,
  there exists at most one inversive distance circle packing metric
  with nonpositive constant combinatorial $R$-curvature.
\end{description}
\end{corollary}

\begin{remark}
The proof of Theorem \ref{theorem global rigidity} follows essentially that given by Luo in \cite{L2} and the case of $\overline{R}=0$ is
in fact included in their result, as $R=0$ is equivalent to $K=0$.
\end{remark}

\begin{remark}\label{uniqueness in R^N}
We can extend the definition of combinatorial curvature to
$\widetilde{R}_i=\frac{\widetilde{K_i}}{s_i^2(r)}$.
If $R(\overline{r})=\overline{R}$ with $\overline{r}\in \Omega$ and $\overline{R}\leq 0$,
under the same conditions as in Theorem \ref{theorem global rigidity},
it can be proved in the same way that
$\overline{r}$ is the unique Euclidean inversive distance circle packing metric in $\mathds{R}^N_{>0}$
with $\widetilde{R}(\overline{r})=\overline{R}$ for the Euclidean
background geometry (up to scaling in the case of $\overline{R}\equiv0$),
and $\overline{r}$ is the unique hyperbolic inversive distance circle packing metric in $\mathds{R}^N_{>0}$
with $\widetilde{R}(\overline{r})=\overline{R}$ for the hyperbolic
background geometry
\end{remark}

\begin{remark}
The condition of $\overline{R}\leq 0$
implies that $\chi(M)\leq 0$ by the combinatorial Gauss-Bonnet formulas (\ref{Gauss-Bonnet with weight}) and (\ref{hyperbolic Gauss-Bonnet formula}).
So Theorem \ref{theorem global rigidity} is in fact only valid for closed triangulated surfaces with nonpositive Euler number.
However, there is no rigidity result for closed triangulated surface with
$\chi(M)>0$, as shown by Example \ref{counter example} in the following,
which implies the rigidity result in Theorem \ref{theorem global rigidity} is sharp in some sense.
\end{remark}

\begin{example} \label{counter example}
Given a topological sphere $\mathbb{S}^2$,
triangulate it into four faces of a single tetrahedron with four vertices $v_i, v_j, v_k, v_l$.
Fix the inversive distance $I\equiv 2$  for every edge in $E$.
Set $r_j=r_k=r_l=x$ and $r_i=1$, then $l_{ij}=l_{ik}=l_{il}=\sqrt{x^2+4x+1}$ and $l_{jk}=l_{kl}=l_{li}=\sqrt{6}x$.
The triangle inequalities are satisfied if and only if $0<x<4+\sqrt{18}$.
By direct calculations, we have $R_i=2\pi-6\arcsin\frac{\sqrt{6}x/2}{\sqrt{x^2+4x+1}}$ and
$R_j=R_k=R_l=\frac{1}{x^2}\left(\frac{2\pi}{3}+2\arcsin \frac{\sqrt{6}x/2}{\sqrt{x^2+4x+1}}\right)$.
Then the triangulated sphere has constant $R$-curvature
if and only if
$$2\pi-6\arcsin\frac{\sqrt{6}x/2}{\sqrt{x^2+4x+1}}=\frac{1}{x^2}\left(\frac{2\pi}{3}+2\arcsin \frac{\sqrt{6}x/2}{\sqrt{x^2+4x+1}}\right),$$
which is equivalent to
\begin{equation}\label{root of function}
\arcsin \frac{\sqrt{6}x/2}{\sqrt{x^2+4x+1}}-\frac{\pi}{3}+\frac{2\pi}{3(3x^2+1)}=0.
\end{equation}
Set the function $f(x)$ to be the term on the left hand side of (\ref{root of function}).
It is obviously that $f(1)=0$, which implies that $r=(1, 1, 1, 1)$ is an inversive distance circle packing metric with
constant $R$-curvature.

Take $x=2$, then $f(2)=\arcsin\sqrt{\frac{6}{13}}-\frac{11\pi}{39}$. Note that $\frac{11\pi}{39}>\frac{\pi}{4}$, we
have $\sin \frac{11\pi}{39}>\frac{\sqrt{2}}{2}>\sqrt{\frac{6}{13}}$, which implies that $f(2)<0$.
Furthermore, $\lim_{x\rightarrow (4+\sqrt{18})^-}f(x)=\frac{\pi}{6}+\frac{2\pi}{9(4+\sqrt{18})^2+3}>0$.
This implies that there exists $x_0\in(2, 4+\sqrt{18})$ such that $f(x_0)=0$ and $r=(1, x_0, x_0, x_0)$ is another inversive distance
circle packing metric which is not scaling equivalent to $(1, 1, 1, 1)$.
In fact, we have the graph for $y=f(x)$ on $[\frac{1}{2}, 8]$ as shown
in figure \ref{Graph of f(x)}, which is plotted by MATLAB.

\begin{figure}[hbp]
\center
\begin{minipage}[t]{0.75\linewidth}
\centering
\includegraphics[width=0.55\textwidth]{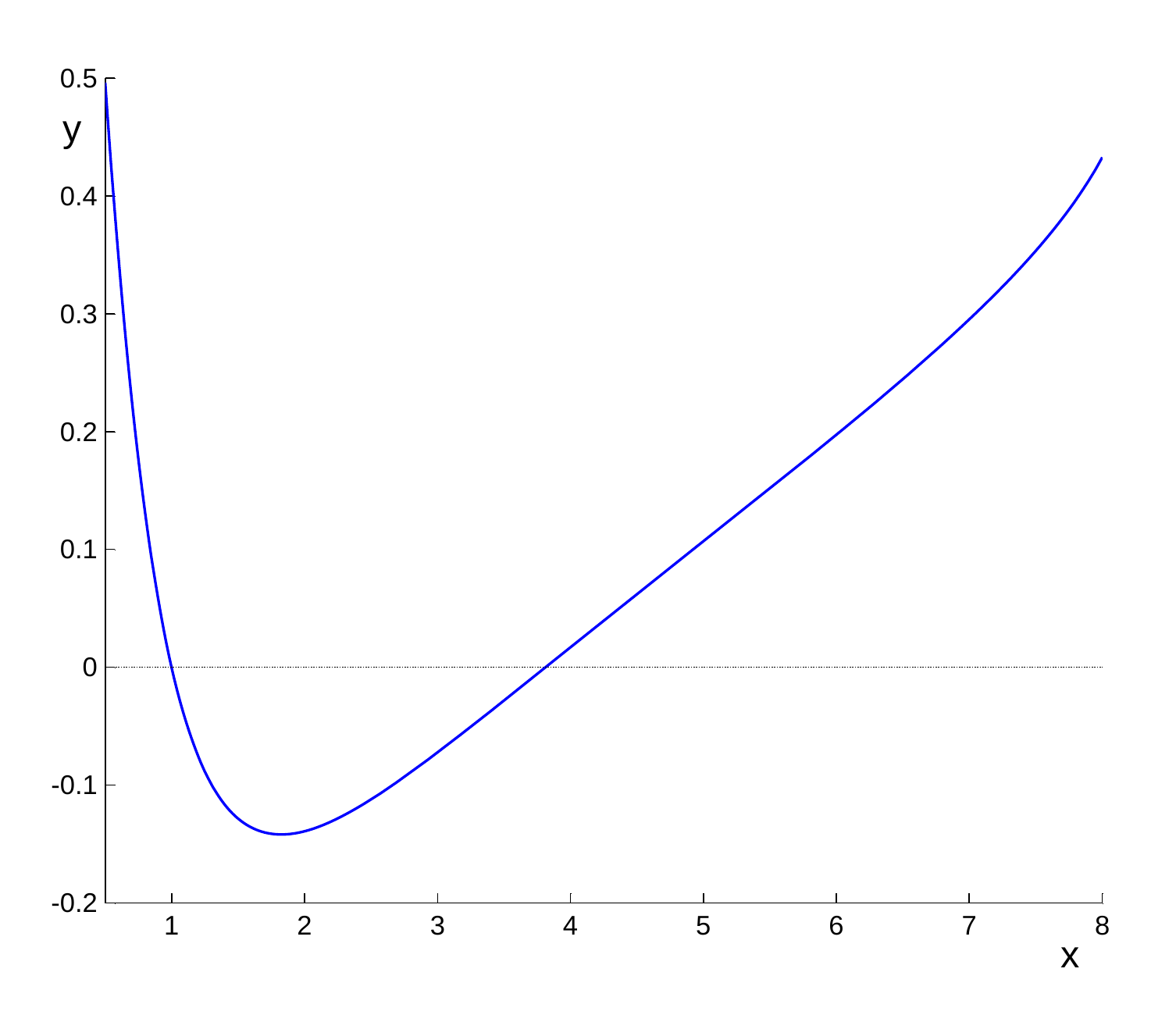}
\caption{{Graph of $f(x)$}}
\label{Graph of f(x)}
\end{minipage}
\end{figure}


So there exist at least two inversive distance circle packing metrics
with the same constant $R$-curvature on $\mathbb{S}^2$ with the tetrahedral triangulation and inversive distance $I\equiv 2$. \qed
\end{example}

\begin{remark}
In \cite{MS}, Ma and Schlenker proved that there is no rigidity for the spherical inversive distance circle packing metrics,
which is also valid for the case we discuss in the paper,
so we focus on the Euclidean and hyperbolic background geometry only.
\end{remark}

\section{Combinatorial Ricci flow}\label{combinatorial Ricci flow}

Given a closed triangulated surface $(M, \mathcal{T})$ with inversive distance $I\geq 0$,
it is natural to consider the following Yamabe problem for the combinatorial curvature $R$.\\

\textbf{Yamabe Problem:} Does there exist an inversive distance circle packing metric with constant combinatorial $R$-curvature or
with prescribed combinatorial $R$-curvature $\overline{R}$? Furthermore, how to find it?\\

Given $\overline{R}\in C(V)$, if there exists some $\overline{r}\in \Omega$ with $R(\overline{r})=\overline{R}$,
we call the function $\overline{R}\in C(V)$ admissible.

Following \cite{GX4}, we introduce an interesting combinatorial curvature flow, the combinatorial Ricci flow, to study the problem.

\subsection{Euclidean combinatorial Ricci flow}\label{subsection Euclidean combinatorial Ricci flow}
As before, we take $r_i^2$ as the analogue of the Riemannian metric for the Euclidean background geometry.
Then we have the following definition of combinatorial Ricci flow.

\begin{definition}
Given a closed triangulated surface $(M, \mathcal{T})$ with inversive distance $I\geq 0$.
For the Euclidean background geometry,
the combinatorial Ricci flow is defined to be
\begin{equation}
\frac{dg_i}{dt}=-R_ig_i
\end{equation}
with normalization
\begin{equation}\label{normalized Euclid Ricci flow}
\frac{dg_i}{dt}=(R_{av}-R_i)g_i,
\end{equation}
where $g_i=r_i^2$ and $R_{av}=\frac{2\pi\chi(M)}{||r||_2^2}$.
\end{definition}

It is easy to check that, along the Euclidean combinatorial Ricci flow (\ref{normalized Euclid Ricci flow}),
the total measure $||r||_2^2$ of the triangulated surface is invariant.
Furthermore, we have the following interesting property for the combinatorial curvature $R_i$
along the normalized Euclidean combinatorial Ricci flow (\ref{normalized Euclid Ricci flow}).

\begin{lemma}\label{evolution of R along Euclidean Ricci flow}
Along the normalized Euclidean combinatorial Ricci flow (\ref{normalized Euclid Ricci flow}), the combinatorial curvature $R$ evolves according to
\begin{equation}\label{evolution equ for R along Ricci flow}
\frac{dR_i}{dt}=\Delta R_i+R_i(R_i-R_{av}),
\end{equation}
where the Laplace operator $\Delta$ is defined to be
\begin{equation}\label{laplace}
\Delta f_i=\frac{1}{r_i^2}\sum_{j=1}^N\left(-\frac{\partial K_i}{\partial u_j}\right)f_j
=\frac{1}{r_i^2}\sum_{j\sim i}\left(-\frac{\partial K_i}{\partial u_j}\right)(f_j-f_i)
\end{equation}
for $f\in C(V)$, where $u_i=\ln r_i^2$.
\end{lemma}
\proof
Note that
$\frac{\partial}{\partial u_j}=\frac{\partial r_j}{\partial u_j}\frac{\partial}{\partial r_j}=\frac{1}{2}r_j\frac{\partial}{\partial r_j}$,
we have
$$\frac{\partial R_i}{\partial u_j}=\frac{\partial}{\partial u_j}(\frac{K_i}{r_i^2})
=\frac{1}{r_i^2}\frac{\partial K_i}{\partial u_j}-\frac{K_i}{r_i^4}\frac{1}{2}r_j\cdot 2r_i\delta_{ij}
=\frac{1}{r_i^2}\frac{\partial K_i}{\partial u_j}-R_i\cdot \frac{r_j}{r_i}\delta_{ij}.$$
Then along the normalized combinatorial Ricci flow (\ref{normalized Euclid Ricci flow}), we have
\begin{equation*}
\begin{aligned}
\frac{dR_i}{dt}=&\sum_j\frac{\partial R_i}{\partial u_j}\frac{d u_j}{d t}\\
=&\sum_j(\frac{1}{r_i^2}\frac{\partial K_i}{\partial u_j}-R_i\cdot \frac{r_j}{r_i}\delta_{ij})(R_{av}-R_j)\\
=&\frac{1}{r_i^2}\sum_j(-\frac{\partial K_i}{\partial u_j})R_j+R_i(R_i-R_{av})\\
=&\Delta R_i+R_i(R_i-R_{av}),
\end{aligned}
\end{equation*}
where Lemma \ref{positivity of L} is used in the third line.
\qed

\begin{remark}
It is very interesting that
the evolution equation (\ref{evolution equ for R along Ricci flow}) for the combinatorial curvature $R$
along the normalized combinatorial Ricci flow (\ref{normalized Euclid Ricci flow})
is a reaction-diffusion equation and
has almost the same form of the evolution of the scalar curvature along the smooth surface Ricci flow
obtained by Hamilton \cite{Ham1}. However, as noted by Guo \cite{Guo}, the coefficient in the Laplacian
operator in (\ref{laplace}) may be negative, which leads to that the discrete maximal principle in \cite{GX4} is not valid here.
\end{remark}

\begin{remark}
The Laplace operator $\Delta$ could be written in the following matrix form
$$\Delta=-\Sigma^{-1}L$$
where $\Sigma=diag\{r_1^2, \cdots, r_N^2\}$. As $\Sigma$ is a diagonal matrix, we have
$$-\Delta=\Sigma^{-1}L\sim \Sigma^{\frac{1}{2}}\cdot \Sigma^{-1}L \cdot \Sigma^{-\frac{1}{2}}=\Sigma^{-\frac{1}{2}}L\Sigma^{-\frac{1}{2}},$$
which is conjugate to $L$ and positive semi-definite. By the property of $L$ in Lemma \ref{positivity of L},
we have $-\Delta$ has one eigenvalue $0$ and all other eigenvalues positive, which is similar
to the case of eigenvalues of Laplacian operators on Riemannian manifolds. So we can also
define the first eigenvalue $\lambda_1(\Delta)$ of the discrete Laplace operator $\Delta$
as the smallest positive number $\lambda$ satisfying $-\Delta f=\lambda f$ for some $f\in C(V)$ with $f\not\equiv0$.
We can also define the $n$th eigenvalue $\lambda_n(\Delta)$ similarly.
\end{remark}

Now we study the long time behavior of the combinatorial Ricci flow. Firstly, we have the following result.

\begin{proposition}\label{necessary for existence of cont curv metric}
Given a closed triangulated surface $(M, \mathcal{T})$ with inversive distance $I\geq 0$.
If the solution of the Euclidean combinatorial Ricci flow (\ref{normalized Euclid Ricci flow}) stays in $\Omega$ and converges
to an inversive distance circle packing metric $r^*\in\Omega$, then there exists a constant $R$-curvature metric in $\Omega$
and $r^*$ is such one.
\end{proposition}
\proof
Suppose $u^*$ is the corresponding $u$-coordinate of $r^*$, i.e. $u^*_i=\ln r^{*2}_i$. Then $R(t)\rightarrow R(u^*)$
by the continuity of the curvature.
As $u(t)\rightarrow u^*\in \mathcal{U}$, there exists $\xi_n\in(n,n+1)$ such that
$$u(n+1)-u(n)=u'(\xi_n)=(R_{av}-R_i)|_{\xi_n}\rightarrow 0$$
as $n\rightarrow +\infty$,
which implies $R(u^*)=R_{av}$ and $r^*$ is a constant curvature metric in $\Omega$.\qed

Proposition \ref{necessary for existence of cont curv metric} shows that the existence of a constant $R$-curvature metric
is necessary for the convergence of the Euclidean combinatorial Ricci flow (\ref{normalized Euclid Ricci flow}).
Conversely, we want to know the behavior of the combinatorial Ricci flow under the existence of a constant curvature metric.
However, the combinatorial Ricci flow may develop singularities, which can be separated into two types as follows.

\begin{definition}\label{definition of singularity}
Given a closed triangulated surface $(M, \mathcal{T})$ with inversive distance $I\geq 0$.
Suppose the solution of the normalized Euclidean combinatorial Ricci flow (\ref{normalized Euclid Ricci flow})
exists on $[0,\tau)$, where $\tau\leq +\infty$.
\begin{description}
  \item[(1)] The normalized Euclidean combinatorial Ricci flow (\ref{normalized Euclid Ricci flow})
  is said to develop an essential singularity at time $\tau$
  if there is a vertex $i\in V$ such that $r_{i}(t_n)$ tends to zero for some time sequence $t_n\rightarrow \tau$.
  \item[(2)] The normalized Euclidean combinatorial Ricci flow (\ref{normalized Euclid Ricci flow}) is said to develop a removable singularity at $\tau$
if there is a sequence of time $t_n\rightarrow \tau$ such that $r(t_n)$ stay in a compact set of $\mathds{R}^N_{>0}$ and
there is a triangle $\triangle ijk\in F$ such that the triangle inequalities degenerate as $t_n\rightarrow \tau$.
\end{description}
\end{definition}

It is interesting that the essential singularity will not develop
along the normalized Euclidean combinatorial Ricci flow (\ref{normalized Euclid Ricci flow}),
under the existence of a nonpositive constant curvature metric. Before presenting the result, we state the following
property for the extended Ricci potential function $\widetilde{F}$.

\begin{proposition}\label{prperty of tilde F}
Given a closed triangulated surface $(M, \mathcal{T})$ with inversive distance $I\geq 0$ and $\chi(M)\leq 0$.
Suppose there exists a Euclidean inversive distance circle packing metric $r^*\in \Omega$ with constant $R$-curvature and
define the extended Ricci potential function $\widetilde{F}(u)$ on $\mathds{R}^N$ as
\begin{equation}
\begin{aligned}
\widetilde{F}(u)=-\sum_{\Delta ijk\in F}\widetilde{F}_{ijk}(u)+\int_{u^*}^u\sum_{i=1}^N(2\pi-R_{av}r_i^2)du_i,
\end{aligned}
\end{equation}
where
$$\widetilde{F}_{ijk}(u)=\int_{u^*}^u\widetilde{\theta}_idu_i+\widetilde{\theta}_jdu_j+\widetilde{\theta}_kdu_k.$$
Then $\widetilde{F}(u+t\textbf{1})=\widetilde{F}(u)$ and $\nabla_{u_i}\widetilde{F}=\widetilde{K}_i-R_{av}r_i^2$,
where $\widetilde{K}_i=2\pi-\sum_{\triangle ijk\in F}\widetilde{\theta}_i$.
Set $\mathcal{U}_a=\{u\in \mathds{R}^N|\sum_{i}u_i=a\}, a\in \mathds{R}$. Then $\widetilde{F}|_{\mathcal{U}_a}$ is convex
and proper. $\widetilde{F}|_{\mathcal{U}_a}$ has a unique zero point which is the unique minimal point of $\widetilde{F}|_{\mathcal{U}_a}$.
Furthermore, $\lim_{u\in \mathcal{U}_a, u\rightarrow \infty}\widetilde{F}(u)=+\infty$.
\end{proposition}
\proof
By the definition of $\tilde{\theta}$, we have $\tilde{\theta}_i+\tilde{\theta}_j+\tilde{\theta}_k=\pi$
for any topological triangle $\triangle ijk\in F$. As the value of $\widetilde{F}_{ijk}(u)$ is independent of the
path from $u^*$ to $u$, it is easy to check $\widetilde{F}_{ijk}(u+t\textbf{1})=\widetilde{F}_{ijk}(u)+t\pi$
by the scaling invariance of $\widetilde{\theta}$ with respect to $r$. Then we have
\begin{equation*}
\begin{aligned}
\widetilde{F}(u+t\textbf{1})-\widetilde{F}(u)
=&-\pi |F|t+\int_{u}^{u+t\textbf{1}}\sum_{i=1}^N(2\pi-R_{av}r_i^2)du_i\\
=&-\pi |F|t+(2\pi |V|-2\pi\chi(M))t\\
=&0.
\end{aligned}
\end{equation*}
It is also easy to check that
$$\nabla_{u_i}\widetilde{F}=-\sum_{\triangle ijk\in F}\widetilde{\theta}_i^{jk}+2\pi-R_{av}r_i^2=\widetilde{K}_i-R_{av}r_i^2$$
and
$$\operatorname{Hess}_{u}\widetilde{F}|_{\mathcal{U}}=L-R_{av}\Sigma^{\frac{1}{2}}\left(I-\frac{rr^T}{||r||^2}\right)\Sigma^{\frac{1}{2}}
=\Sigma^{\frac{1}{2}}\left(\Lambda-R_{av}\left(I-\frac{rr^T}{||r||^2}\right)\right)\Sigma^{\frac{1}{2}},$$
where $\Sigma=diag\{r_1^2, \cdots, r_N^2\}$,
$\Lambda=\Sigma^{-\frac{1}{2}}L\Sigma^{-\frac{1}{2}}$.
Note that $\Lambda$ is positive semi-definite with $rank=N-1$ and kernel $\{cr|c\in \mathds{R}\}$ by Lemma \ref{positivity of L}, we have
$\operatorname{Hess}_{u}\widetilde{F}|_{\mathcal{U}}$ is positive semi-definite with $rank=N-1$ and kernel $\{c\textbf{1}|c\in \mathds{R}\}$
under the condition $\chi(M)\leq 0$. This implies that $\operatorname{Hess}_{u}\widetilde{F}|_{\mathcal{U}\cap \mathcal{U}_a}$ is
strictly positive definite and $\widetilde{F}|_{\mathcal{U}\cap \mathcal{U}_a}$ is strictly convex.
By the concavity of $\widetilde{F}_{ijk}$ and the condition $\chi(M)\leq 0$, it is easy to check that $\widetilde{F}$ is $C^1$ smooth
and convex. Combining with $\widetilde{F}(u+t\textbf{1})=\widetilde{F}(u)$, we have
$\widetilde{F}|_{\mathcal{U}_a}$ is a $C^1$ smooth convex function and strictly convex on $\mathcal{U}\cap \mathcal{U}_a$
under the condition $\chi(M)\leq 0$.

The existence of a constant curvature metric $r^*$ ensures that the corresponding $u^*\in \mathcal{U}\cap \mathcal{U}_a$ is
a critical point of $\widetilde{F}|_{\mathcal{U}_a}$.
Then the rest of the proof follows from applying the following Lemma \ref{lemma on convex function}
to the convex function $\widetilde{F}|_{\mathcal{U}_a}$.
%
%
%
\qed

\begin{lemma}\label{lemma on convex function}
Suppose $f(x)$ is a $C^1$ smooth convex function on $\mathds{R}^n$ with $\nabla f(x_0)=0$ for some $x_0\in \mathds{R}^n$,
$f(x)$ is $C^2$ smooth and strictly convex
in a neighborhood of $x_0$, then $\lim_{x\rightarrow \infty}f(x)=+\infty$.
\end{lemma}
\proof
Without loss of generality, we can assume $x_0=0$, otherwise parallel translating $x_0$ to $0$.
Suppose $f(x)$ is $C^2$ smooth and strictly convex on the ball $B(0, \epsilon)$ for some constant $\epsilon>0$.
For any $v\in \mathbb{S}^{n-1}$, we have $f_{v}(t)=f(tv)$ is a $C^1$ smooth convex function in $t$ and $f_v'(0)=0$.
Furthermore, $f_{v}(t)$ is strictly convex on $[0, \epsilon]$ since $f(x)$ is strictly convex on $B(0, \epsilon)$.
Then, for $t>\epsilon$, we have
$$f_{v}(t)\geq f_{v}(\epsilon)+f_{v}'(\epsilon)(t-\epsilon).$$
Therefore, for any $x\in \mathds{R}^n$ with $||x||>\epsilon$, we have
\begin{equation}\label{a}
f(x)=f(||x||\frac{x}{||x||})=f_{\frac{x}{||x||}}(||x||)\geq f_{\frac{x}{||x||}}(\epsilon)+f_{\frac{x}{||x||}}'(\epsilon)(||x||-\epsilon).
\end{equation}
Note that
\begin{equation}\label{b}
f_{\frac{x}{||x||}}(\epsilon)=f(\epsilon \frac{x}{||x||})\geq min_{v\in \mathbb{S}^{n-1}}f(\epsilon v)=c_1>0.
\end{equation}
Furthermore, for any $v\in \mathbb{S}^{n-1}$, $f_{v}'(\epsilon)>0$ by the strict convexity of $f_v(t)$ on $[0, \epsilon]$ and $f_v'(0)=0$,
so we have
\begin{equation}\label{c}
f_{\frac{x}{||x||}}'(\epsilon)\geq min_{v\in \mathbb{S}^{n-1}}f_{v}'(\epsilon)=c_2>0.
\end{equation}
Combining (\ref{a}), (\ref{b}) and (\ref{c}), we have
$$f(x)\geq c_1+c_2(||x||-\epsilon)$$
for some $c_1, c_2>0$ and any $x\in \mathds{R}^n$ with $||x||>\epsilon$,
which implies $\lim_{x\rightarrow \infty}f(x)=+\infty$. \qed

\begin{remark}\label{uniqueness of const curv metric in R^N}
In fact, the proof of Proposition \ref{prperty of tilde F} shows that, under the conditions of Proposition \ref{prperty of tilde F},
the constant curvature metric $r^*\in \Omega$ is the unique inversive distance circle packing metric with constant $\widetilde{R}$-curvature in $\mathds{R}^N_{>0}$
(up to scaling in the case of $\chi(M)=0$).
\end{remark}


Note that $||r||_2^2=\sum_{i=1}^2$ is invariant along the (\ref{normalized Euclid Ricci flow}), we assume that
$r(0)\in \mathbb{S}^{N-1}\cap \Omega$ in the following.
Set
$$\mathcal{S}=\{u\in \mathds{R}^N|\sum_{i=1}^Ne^{u_i}=1\}.$$
Using the property $\widetilde{F}(u+t\textbf{1})=\widetilde{F}(u)$,
we have the following property for the extended Ricci potential $\widetilde{F}(u)$.

\begin{corollary}
Given a closed triangulated surface $(M, \mathcal{T})$ with inversive distance $I\geq 0$ and $\chi(M)\leq 0$.
Suppose there exists a Euclidean inversive distance circle packing metric $r^*\in \Omega$ with constant $R$-curvature. Then
$\widetilde{F}|_{\mathcal{S}}$ is proper. Furthermore, $\lim_{u\rightarrow \infty}\widetilde{F}|_{\mathcal{S}}=+\infty$.
\end{corollary}
\proof
Let $\Pi$ be the orthogonal projection from $\mathcal{S}$ to the plane $\mathcal{U}_0$.
Then for a sequence $u^{(n)}\in \mathcal{S}$, $\Pi(u^{(n)})$ is unbounded
if and only if $|\Pi(u^{(n)})_i-\Pi(u^{(n)})_j|=|u^{(n)}_i-u^{(n)}_j|$
is unbounded.
If $u\in\mathcal{S}$ and $u\rightarrow \infty$, we have $\Pi(u)\rightarrow \infty$.
Then $\widetilde{F}(u)=\widetilde{F}(\Pi(u))\rightarrow \infty$.
\qed

Now we can get the following result on nonexistence of essential singularities
along the normalized Euclidean combinatorial Ricci flow.

\begin{theorem}\label{nonexistence of essential singularity}
Given a closed triangulated surface $(M, \mathcal{T})$ with inversive distance $I\geq 0$ and $\chi(M)\leq 0$.
If there exists a Euclidean inversive distance circle packing metric $r^*\in \Omega$ with constant $R$-curvature,
then the solution of the combinatorial
Ricci flow (\ref{normalized Euclid Ricci flow}) develop no essential singularities.
\end{theorem}
\proof
To prove the result, we need to extend the combinatorial Ricci flow to the following form
\begin{equation}\label{extended ricci flow}
\frac{dg_i}{dt}=(R_{av}-\widetilde{R}_i)g_i,
\end{equation}
where $\widetilde{R}_i=\frac{\widetilde{K}_i}{r_i^2}$ with $\widetilde{K}_i=2\pi-\sum_{\Delta ijk\in F}\widetilde{\theta}_i^{jk}$.
The equation (\ref{extended ricci flow}) could be written as
\begin{equation}\label{extended ricci flow in u}
\frac{du_i}{dt}=R_{av}-\widetilde{R}_i.
\end{equation}
Note that $R_{av}-\widetilde{R}_i$ is just a continuous function of $u\in \mathds{R}^N$ and not locally
Lipschiz in $u\in \mathds{R}^N$. By the existence of solution of ODE for continuous function, there exists at least one
solution of (\ref{extended ricci flow}). There may be no uniqueness for the solution of the extended flow (\ref{extended ricci flow})
on $\mathds{R}^N$. However, note that $\widetilde{R}|_{\mathcal{U}}=R$ on $\mathcal{U}$ and $R$ is locally Lipschitz,
so the solutions of (\ref{normalized Euclid Ricci flow}) and (\ref{extended ricci flow in u}) agree on $\mathcal{U}$ by
Picard's uniqueness for the solution of ODE. And then any solution of (\ref{extended ricci flow in u})
extends the solution of (\ref{normalized Euclid Ricci flow}).

Under the assumption of the existence of the constant curvature metric,
the Ricci potential
$$
\widetilde{F}(u)=-\sum_{\Delta ijk\in F}\widetilde{F}_{ijk}(u)+\int_{u^*}^u\sum_{i=1}^N(2\pi-R_{av}r_i^2)du_i,
$$
is convex.
Note that $\widetilde{F}|_{\mathcal{S}}$ is proper, $\lim_{u\rightarrow \infty}\widetilde{F}|_{\mathcal{S}}(u)=+\infty$ and
$$\frac{d}{dt}\widetilde{F}(u(t))=-\sum_{i=1}^N(\widetilde{R}_i-R_{av})^2r_i^2\leq 0$$
along the extended Ricci flow (\ref{extended ricci flow in u}),
we have the solution $u(t)$ of the extended flow (\ref{extended ricci flow in u}) stays in a compact subset of $\mathds{R}^N$, which
implies that the solution $r_i(t)$ of (\ref{extended ricci flow}) has uniform positive lower and upper bounds for
any $i\in V$ and $t\in \mathds{R}$.
As the solution of (\ref{extended ricci flow}) extends the solution of the
combinatorial Ricci flow (\ref{normalized Euclid Ricci flow}), we have
the solution of the combinatorial Ricci flow (\ref{normalized Euclid Ricci flow}) will not develop essential singularities.\qed
%

Theorem \ref{nonexistence of essential singularity} shows that there is no essential singularities along the
normalized Euclidean combinatorial Ricci flow (\ref{normalized Euclid Ricci flow}) in finite time and at time infinity.
If the flow (\ref{normalized Euclid Ricci flow}) further develops no removable singularities, we have the following interesting result.

\begin{theorem}\label{convergence under nonexistence of removable singularity}
Given a closed triangulated surface $(M, \mathcal{T})$ with inversive distance $I\geq 0$ and $\chi(M)\leq 0$.
If there exists a Euclidean inversive distance circle packing metric $r^*\in \Omega$ with constant $R$-curvature
and the solution of the combinatorial
Ricci flow (\ref{normalized Euclid Ricci flow}) develops no removable singularities in finite time,
then the solution converges exponentially fast to the constant curvature metric.
\end{theorem}
\proof
The proof of Theorem (\ref{nonexistence of essential singularity}) shows that the solution of
(\ref{normalized Euclid Ricci flow}) and (\ref{extended ricci flow}) stays in a compact subset of $\mathds{R}^N_{>0}$.
This implies the solution of (\ref{normalized Euclid Ricci flow}) and (\ref{extended ricci flow})
exists for all time.
$u(t)$ is bounded implies that the function $\widetilde{F}(u(t))$ is bounded along the flow.
Note that
$$\frac{d}{dt}\widetilde{F}(u(t))=\sum_{i=1}^N\nabla_{u_i}\widetilde{F}\cdot \frac{du_i}{dt}
=\sum_{i=1}^N (K_i-R_{av}r_i^2)(R_{av}-R_i)=-\sum_{i=1}^N (R_{av}-R_i)^2r_i^2\leq 0,$$
we have $\lim_{t\rightarrow +\infty}\widetilde{F}(u(t))$ exists. Then
\begin{equation}\label{derivative of F(u(t))}
\begin{aligned}
\widetilde{F}(u(n+1))-\widetilde{F}(u(n))=-\sum_{i=1}^N (R_{av}-\widetilde{R}_i)^2r_i^2|_{t=\xi_n}\rightarrow 0
\end{aligned}
\end{equation}
as $n\rightarrow +\infty$, for some $\xi_n\in (n, n+1)$. As $r(\xi_n)$ is bounded, then there exists a subsequence,
still denoted as $r(\xi_n)$, such that $r(\xi_n)\rightarrow \bar{u}$ and then $\widetilde{R}(r(\xi_n))\rightarrow \widetilde{R}(\bar{r})$.
The equation (\ref{derivative of F(u(t))}) shows that $\widetilde{R}(\bar{r})= R_{av}$.
By the uniqueness of constant curvature metric in $\mathds{R}^N_{>0}$ in Remark \ref{uniqueness of const curv metric in R^N},
we have $\bar{r}=r^*\in \Omega$.

Furthermore, $r^*$ is a local attractor of the equation  (\ref{normalized Euclid Ricci flow}).
In fact, the equation (\ref{normalized Euclid Ricci flow}) could be written as
$$\frac{dr_i}{dt}=\frac{1}{2}(R_{av}-R_i)r_i.$$
Set $\Gamma_{i}(r)=\frac{1}{2}(R_{av}-R_i)r_i$, then we have
\begin{equation*}
D_r\Gamma|_{r^*}=R_{av}\left(I-\frac{r r^T}{||r||^2}\right)-\Sigma^{-\frac{1}{2}} L \Sigma^{-\frac{1}{2}}=R_{av}\left(I-\frac{r r^T}{||r||^2}\right)-\Lambda,
\end{equation*}
where $\Sigma=diag\{r_1^2,\cdots,r_N^2\}$ and
$\Lambda=\Sigma^{-\frac{1}{2}} L \Sigma^{-\frac{1}{2}}\sim \Sigma^{-\frac{1}{2}}\Lambda\Sigma^{\frac{1}{2}}=-\Delta$.
Select an orthonormal matrix $P$ such that
$$P^T\Lambda P=diag\{0,\lambda_1(\Lambda),\cdots,\lambda_{N-1}(\Lambda)\}.$$
Suppose $P=(e_0,e_1,\cdots,e_{N-1})$,
where $e_i$ is the $(i+1)$-column of $P$.
Then $\Lambda e_0=0$ and $\Lambda e_i=\lambda_i e_i,\,1\leq i\leq N-1$,
which implies $e_0=r/\|r\|$ and $r\perp e_i,\,1\leq i\leq N-1$.
Hence $\big(I_{N}-\frac{rr^T}{\|r\|^{2}}\big)e_0=0$ and $\big(I_{N}-\frac{rr^T}{\|r\|^{2}}\big)e_i=e_i$, $1\leq i\leq N-1$,
which implies $P^T\big(I_{N}-\frac{rr^T}{\|r\|^{2}}\big)P=diag\{0,1,\cdots,1\}$. Therefore,
\begin{equation*}
D_{r}\Gamma\big|_{r^*}=
P \cdot diag\{0,R_{av}-\lambda_1(\Lambda),\cdots,R_{av}-\lambda_{N-1}(\Lambda)\}\cdot P^T.
\end{equation*}
Note that $\lambda_1(\Lambda)>R_{av}^*=\frac{2\pi\chi(M)}{||r^*||^2}$ in the case of $\chi(M)\leq 0$,
then $D_{r}\Gamma\big|_{r^*}$ is negative semi-definite with kernel $\{cr| c\in \mathbb{R}\}$ and $rank\,(D_{r}\Gamma\big|_{r^*})=N-1$.
Note that, along the flow (\ref{normalized Euclid Ricci flow}), $||r||^2$ is invariant. Thus the kernel is transversal to the flow. This implies that $D_r\Gamma|_{r^*}$ is negative definite on $\mathbb{S}^{N-1}$ and $r^*$ is a local attractor of the normalized combinatorial Ricci
flow (\ref{normalized Euclid Ricci flow}).
Then the conclusion follows from the Lyapunov Stability Theorem(\cite{P}, Chapter 5) and the fact that $r(\xi_n)\rightarrow r^*$.
\qed

\begin{remark}
The proof of Theorem \ref{convergence under nonexistence of removable singularity} shows that
the solution of (\ref{normalized Euclid Ricci flow}) will not develop any removable singularity at time infinity.
\end{remark}

\begin{remark}
Note that $\Lambda=\Sigma^{-\frac{1}{2}} L \Sigma^{-\frac{1}{2}}\sim \Sigma^{-\frac{1}{2}}\Lambda\Sigma^{\frac{1}{2}}=-\Delta$,
to ensure the local convergence of the flow (\ref{normalized Euclid Ricci flow}), we only need the first eigenvalue
of the discrete Laplace operator satisfying $\lambda_1(\Delta)|_{r^*}>R_{av}$, which is obvious in the case of $\chi(M)\leq 0$.
\end{remark}

The only case left is that the solution of the normalized combinatorial Ricci
flow (\ref{normalized Euclid Ricci flow}) develops removable singularities at finite time,
which implies that there exists some triangle $\triangle ijk\in F$ degenerating as $t\rightarrow \tau$.
The proof of Theorem \ref{convergence under nonexistence of removable singularity}
suggests us to use the extended combinatorial Ricci flow (\ref{extended ricci flow})
to characterize the existence of the constant curvature metric.
In fact, we have the following main theorem.

\begin{theorem}\label{convergence of extended flow}
Given a closed triangulated surface $(M, \mathcal{T})$ with inversive distance $I\geq 0$ and $\chi(M)\leq 0$.
If there exists an inversive distance circle packing metric $r^*\in \Omega$ with constant $R$-curvature,
then the solution of the extended combinatorial Ricci flow (\ref{extended ricci flow}) exists for all time and converges
exponentially fast to the constant curvature metric $r^*$.
\end{theorem}
The proof of Theorem \ref{convergence of extended flow} is almost the same
as that of Theorem \ref{convergence under nonexistence of removable singularity}, we omit the details here.
And in the special case that the solution of the normalized combinatorial Ricci flow (\ref{normalized Euclid Ricci flow})
stays in $\Omega$, this is just Theorem \ref{convergence under nonexistence of removable singularity}.

\begin{remark}
As there maybe exist more than one solution for the extended flow (\ref{extended ricci flow}) and
every solution of (\ref{extended ricci flow}) extends the solution of the combinatorial Ricci flow (\ref{normalized Euclid Ricci flow}),
Theorem  \ref{convergence of extended flow} shows that every such solution will go back into $\Omega$ and converges
to the constant curvature metric. Furthermore, as $\partial\Omega$ is piecewise analytic and the
solution of the extended flow (\ref{extended ricci flow}) is $C^1$ smooth, any solution of (\ref{extended ricci flow})
intersects $\partial \Omega$ at most finite times.
\end{remark}

\begin{remark}
The existence of removable singularity along the normalized combinatorial Ricci flow
suggests a way to do surgery on the triangulated surface, as stated by Luo in \cite{L1}.
And this way should be also valid here.
However, surgeries on the triangulated surface change the combinatorial structure on the
triangulated surface. If there exists a constant curvature metric on the triangulated surface
with given combinatorial structure and the solution of the extended flow (\ref{extended ricci flow})
develops a removable singularity at finite time, surgeries on the triangulated surface shall
give a constant curvature metric with a different combinatorial structure.
However, Theorem  \ref{convergence of extended flow} ensures that we can find the constant curvature metric $r^*\in \Omega$
with the given combinatorial structure under the existence of a constant curvature metric, even if the
solution of the normalized combinatorial Ricci flow (\ref{normalized Euclid Ricci flow}) develops
removable singularities at finite time.
\end{remark}

We can also use the Euclidean combinatorial flow to study the prescribing curvature problem of the $R$-curvature.
As the proof is almost parallel to the case of constant curvature metric problem, we just
summarize the results in the following theorem and omit its proof.

\begin{theorem}\label{Euclidean prescribing curvature problem}
Given a closed triangulated surface $(M, \mathcal{T})$ with inversive distance $I\geq 0$ and $\chi(M)\leq 0$.
$\overline{R}\in C(V)$ is a given function on $(M, \mathcal{T})$.
\begin{description}
  \item[(1)] If the solution of the modified Euclidean combinatorial Ricci flow
  \begin{equation}\label{modified Euclidean combinatorial Ricci flow}
  \frac{dg_i}{dt}=(\overline{R}_i-R_i)g_i,
  \end{equation}
where $g_i=r_i^2$,
  stays in $\Omega$ and converges to an inversive distance circle packing metric $\overline{r}\in \Omega$,
  then $R(\overline{r})=\overline{R}$ and
  $\overline{R}$ is admissible.
  \item[(2)] Suppose $\overline{R}\in C(V)$ is admissible with $R(\overline{r})=\overline{R}$ for some $\overline{r}\in \Omega$,
$\overline{R}\leq0$ and $\overline{R}\not\equiv0$,
  then the modified Euclidean combinatorial Ricci flow (\ref{modified Euclidean combinatorial Ricci flow})
  does not develop essential singularities in finite time and at
  time infinity. Furthermore, if the solution of (\ref{modified Euclidean combinatorial Ricci flow})
  develops no removable singularities in finite time, then the solution of (\ref{modified Euclidean combinatorial Ricci flow})
  exists for all time,
  converges exponentially fast to $\overline{r}$ and does not develop removable singularities at time infinity.
  \item[(3)] Suppose $\overline{R}\in C(V)$ is admissible with $R(\overline{r})=\overline{R}$ for some $\overline{r}\in \Omega$,
$\overline{R}\leq0$ and $\overline{R}\not\equiv0$,
  then the solution of the extended modified
  Euclidean combinatorial Ricci flow
  $$\frac{dg_i}{dt}=(\overline{R}_i-\widetilde{R}_i)g_i$$
  exists for all time and converges exponentially fast to $\overline{r}$.
\end{description}
\end{theorem}

\begin{remark}
The case of $\overline{R}\equiv 0$ is just the constant curvature problem,
so we do not include it in Theorem \ref{Euclidean prescribing curvature problem} here.
Theorem \ref{Euclidean prescribing curvature problem} provided an effective way to find an inversive distance circle
packing metric with given prescribed curvature and given combinatorial structure.
\end{remark}

\subsection{Hyperbolic combinatorial Ricci flow}\label{subsection hyperbolic combinatorial Ricci flow}
For the hyperbolic background geometry, we take $\tanh^2\frac{r_i}{2}$
as the analogue of the Riemannian metric. Then we have the following definition of
hyperbolic combinatorial Ricci flow.

\begin{definition}
Given a closed triangulated surface $(M, \mathcal{T})$ with inversive distance $I\geq 0$.
For the hyperbolic background geometry, the combinatorial Ricci flow is defined as
\begin{equation}
\frac{dg_i}{dt}=-R_ig_i,
\end{equation}
where $g_i=\tanh^2\frac{r_i}{2}$.
Given a function $\overline{R}\in C(V)$, the modified combinatorial Ricci flow is defined to be
\begin{equation}\label{modifed hyperbolic Ricci flow}
\frac{dg_i}{dt}=(\overline{R}_i-R_i)g_i.
\end{equation}
\end{definition}

Similar to the Euclidean background geometry,
we have the following evolution of the hyperbolic combinatorial curvature $R$
along the modified hyperbolic combinatorial Ricci flow.

\begin{lemma}
Along the modified hyperbolic combinatorial Ricci flow (\ref{modifed hyperbolic Ricci flow}), the combinatorial curvature $R_i$ evolves according to
\begin{equation*}
\begin{aligned}
\frac{dR_i}{dt}=-\frac{1}{\tanh^2\frac{r_i}{2}}\sum_{j}\frac{\partial K_i}{\partial u_j}(R_j-\overline{R}_j)+R_i(R_i-\overline{R}_i),
\end{aligned}
\end{equation*}
where $u_i=\ln \tanh^2\frac{r_i}{2}$.
\end{lemma}
\proof
As $u_i=\ln \tanh^2\frac{r_i}{2}$, we have
$$du_i=\frac{2}{\tanh\frac{r_i}{2}}\cdot \frac{1}{\cosh^2\frac{r_i}{2}}\cdot \frac{1}{2}dr_i=\frac{2}{\sinh r_i}dr_i,$$
which implies that
$\frac{\partial}{\partial u_j}=\frac{\partial r_j}{\partial u_j}\frac{\partial}{\partial r_j}=\frac{1}{2}\sinh r_j\frac{\partial}{\partial r_j}$.
So we have
\begin{equation*}
\begin{aligned}
\frac{dR_i}{dt}
=&\sum_j\frac{\partial R_i}{\partial u_j}\frac{du_j}{dt}\\
=&\sum_j\left(\frac{1}{\tanh^2\frac{r_i}{2}}\frac{\partial K_i}{\partial u_j}-\frac{1}{2}\sinh r_j\cdot\frac{R_i}{\tanh^2\frac{r_i}{2}}\cdot2\tanh\frac{r_i}{2}\frac{1}{\cosh^2\frac{r_i}{2}}\cdot\frac{1}{2}\delta_{ij}\right)(\overline{R}_j-R_j)\\
=&\sum_j\left(\frac{1}{\tanh^2\frac{r_i}{2}}\frac{\partial K_i}{\partial u_j}-\frac{\sinh r_j}{\sinh r_i}R_i\delta_{ij}\right)(\overline{R}_j-R_j)\\
=&-\frac{1}{\tanh^2\frac{r_i}{2}}\sum_{j}\frac{\partial K_i}{\partial u_j}(R_j-\overline{R}_j)+R_i(R_i-\overline{R}_i).
\end{aligned}
\end{equation*}
\qed

Different from the case of Euclidean background geometry, we just consider the prescribing curvature problem in the case of
hyperbolic background geometry using the modified hyperbolic combinatorial Ricci flow (\ref{modifed hyperbolic Ricci flow})
following \cite{GX3}, which includes the constant curvature problem with the constant being zero.
First, we have the following necessary condition for the convergence of the flow.

\begin{proposition}\label{necessary for the convergence of the midified hyper flow}
Given a closed triangulated surface $(M, \mathcal{T})$ with inversive distance $I\geq 0$.
If the solution of the modified hyperbolic combinatorial Ricci flow (\ref{modifed hyperbolic Ricci flow})
stays in $\Omega$ and converges to $\overline{r}\in \Omega$, then we have $\overline{R}$ is admissible and
$R(\overline{r})=\overline{R}$.
\end{proposition}

The proof of Proposition \ref{necessary for the convergence of the midified hyper flow}
is similar to that of Proposition \ref{necessary for existence of cont curv metric}, we omit it here.

Conversely, assuming the admissibility of $\overline{R}$, we have results similar to that of the Euclidean background case.
Before presenting the results, we state the following useful lemma which we need in the following,
which is given in \cite{CL1, GX3} and still valid in the case of
inversive distance circle packing metric with $I\geq 0$.

\begin{lemma}\label{limit behavior of theta_i}
Given a hyperbolic triangle $\triangle v_iv_jv_k$, which is constructed by inversive distance circle packing metric $r$ with
inversive distance $I\geq 0$. For any $\epsilon>0$, there exists a number $l>0$ such that when $r_i>l$, the inner angle
$\theta_{i}$ at the vertex $v_i$ of the hyperbolic triangle $\triangle v_iv_jv_k$ is smaller than $\epsilon$.
\end{lemma}

The proof of Lemma \ref{limit behavior of theta_i} is the same as that of Lemma 3.2 in \cite{GX3}, we omit it here.
Note that the convergence result in  Lemma \ref{limit behavior of theta_i} is in fact
the uniform convergence of $\theta_i^{jk}$ with respect to $r_i$, assuming that the triangle inequalities are satisfied.

We further have the following interesting result on triangle inequalities.

\begin{lemma}\label{triangle inequality}
Given a topological triangle $\triangle ijk\in F$ with the lengths
given by (\ref{definition of length of edge for hyperbolic background}) and inversive distance $I\geq 0$,
there exists a positive number $l=l(I_{jk})$ such that if $r_i>l$, then $l_{ij}+l_{ik}>l_{jk}$.
\end{lemma}
\proof
We just need to prove $\cosh(l_{ij}+l_{ik})>\cosh l_{jk}$.
Note that there exists $l>0$ such that if $r_i>l$, then $\cosh^2r_i>I_{jk}+1$. Then we have
\begin{equation*}
\begin{aligned}
\cosh(l_{ij}+l_{ik})
>&\cosh l_{ij}\cosh l_{ik}\\
=&(\cosh r_i\cosh r_j+I_{ij}\sinh r_i\sinh r_j)(\cosh r_i\cosh r_k+I_{ik}\sinh r_i\sinh r_k)\\
>&\cosh^2r_i\cosh r_j\cosh r_k\\
>&(I_{jk}+1)\cosh r_j\cosh r_k\\
\geq& \cosh r_j\cosh r_k+I_{jk}\sinh r_j\sinh r_k\\
=&\cosh l_{jk}.
\end{aligned}
\end{equation*}
\qed

Combining Lemma \ref{limit behavior of theta_i} and Lemma \ref{triangle inequality}, we have

\begin{lemma}\label{limit behavior of extended theta}
Given a topological triangle $\triangle ijk\in F$ with the lengths
given by (\ref{definition of length of edge for hyperbolic background}) with inversive distance $I\geq 0$.
For any $\epsilon>0$,
there exists a positive number $l=l(I_{jk}, \epsilon)$ such that if $r_i>l$,
then the extended inner angle $\widetilde{\theta}_i^{jk}<\epsilon$.
\end{lemma}

\proof
By Lemma \ref{triangle inequality}, there exists $l_1>0$ such that $l_{ij}+l_{ik}>l_{jk}$.
If the other two triangle inequalities are not satisfied, then $\widetilde{\theta}_i^{jk}=0$.
If the other two triangle inequalities are satisfied, by Lemma \ref{limit behavior of theta_i},
there exists $l_2>0$ such that if $r_i>l_2$, then $\widetilde{\theta}_i^{jk}<\epsilon$.
Take $l=\max\{l_1, l_2\}$, then we get the proof of the lemma.\qed

Then we have the following main result.

\begin{theorem}\label{hyperbolic prescribing curvature problem}
Given a closed triangulated surface $(M, \mathcal{T})$ with inversive distance $I\geq 0$ and $\chi(M)\leq 0$.
$\overline{R}\in C(V)$ is admissible with $R(\overline{r})=\overline{R}$ and $\overline{R}\leq 0$.
\begin{description}
  \item[(1)] The modified hyperbolic combinatorial Ricci flow (\ref{modifed hyperbolic Ricci flow})
  does not develop essential singularity in finite time and at time infinity.
  \item[(2)] If the solution of (\ref{modifed hyperbolic Ricci flow})
  develops no removable singularities in finite time, then the solution of (\ref{modifed hyperbolic Ricci flow})
  exists for all time,
  converges exponentially fast to $\overline{r}$ and does not develop removable singularities at time infinity.
  \item[(3)] Any solution of the extended hyperbolic flow
  $$\frac{dg_i}{dt}=(\overline{R}_i-\widetilde{R}_i)g_i$$
  with $g_i=\tanh^2\frac{r_i}{2}$ exists for all time and converges exponentially fast to $\overline{r}$.
\end{description}
\end{theorem}
\proof
\begin{description}
  \item[(1)]
Set $u_i=\ln \tanh^2\frac{r_i}{2}$, which maps $\mathds{R}^N_{>0}$ onto $\mathds{R}^N_{<0}$.
Denote $\overline{u}$ as the $u$-coordinate of $\overline{r}$.
Consider the Ricci potential function
$$F(u)=-\sum_{\triangle ijk\in F}F_{ijk}(u)+\int_{\overline{u}}^u\sum_{i=1}^N(2\pi-\overline{R}_i\tanh^2\frac{r_i}{2}),$$
which is strictly convex on $\mathcal{U}$ in the case of $\overline{R}\leq 0$ and
could be extended to a $C^1$ smooth convex function
$$\widetilde{F}(u)=-\sum_{\triangle ijk\in F}\widetilde{F}_{ijk}(u)+
\int_{\overline{u}}^u\sum_{i=1}^N(2\pi-\overline{R}_i\tanh^2\frac{r_i}{2})$$
defined on $\mathds{R}^N_{<0}$. By direct calculations, we have
$$\nabla_{u_i}\widetilde{F}=\widetilde{K}_i-\overline{R}_i\tanh^2\frac{r_i}{2}$$
and $\overline{u}$ is a critical point of $\widetilde{F}$.
Note that, in the case of $\overline{R}\leq 0$, we have
$$\operatorname{Hess}_u\widetilde{F}|_{\mathcal{U}}=L-\left(
                                                        \begin{array}{ccc}
                                                          \overline{R}_1\tanh^2\frac{r_1}{2} &   &   \\
                                                            & \ddots &   \\
                                                            &   & \overline{R}_N\tanh^2\frac{r_N}{2} \\
                                                        \end{array}
                                                      \right),
$$
which is positive definite by Lemma \ref{positivity of L}. Applying Lemma \ref{lemma on convex function},
we have $\lim_{u\rightarrow \infty}\widetilde{F}(u)=+\infty$.
Along the extended modified hyperbolic combinatorial Ricci flow
\begin{equation}\label{extended modified hyperbolic combinatorial Ricci flow}
\frac{dg_i}{dt}=(\overline{R}_i-\widetilde{R}_i)g_i,
\end{equation}
where $g_i=\tanh^2\frac{r_i}{2}$,
we have
\begin{equation*}
\frac{d}{dt}\widetilde{F}(u(t))=\nabla_u\widetilde{F}\cdot \frac{du}{dt}=-\sum_{i=1}^N(\overline{R}_i-R_i)^2\tanh^2\frac{r_i}{2}\leq 0,
\end{equation*}
which implies $\widetilde{F}(u(t))$ is decreasing along the extended hyperbolic flow (\ref{extended modified hyperbolic combinatorial Ricci flow}).
By $\lim_{u\rightarrow \infty}\widetilde{F}(u)=+\infty$, we have the solution $u(t)$ of the
extended hyperbolic flow (\ref{extended modified hyperbolic combinatorial Ricci flow})
is bounded in $\mathds{R}^N_{<0}$, and then there exists $c>0$ such that $u_i(t)>-c$ for any $i\in V$,
which implies $r_i(t)\geq \delta=\ln \frac{1+e^{-c}}{1-e^{-c}}>0$.
As the solution of the hyperbolic flow (\ref{modifed hyperbolic Ricci flow}) is always part of
the solution of the extended hyperbolic flow (\ref{extended modified hyperbolic combinatorial Ricci flow}),
we have the solution of the hyperbolic flow (\ref{modifed hyperbolic Ricci flow})
never develops essential singularities.
  \item[(2)]
In the case that the solution of the hyperbolic Ricci flow (\ref{modifed hyperbolic Ricci flow})
develops no removable singularities at finite time,
the solution of the hyperbolic Ricci flow (\ref{modifed hyperbolic Ricci flow}) exists for all time.
We claim that $r_i(t)$ is bounded from above.
Otherwise there exists at least one vertex $i\in V$ such that $\limsup_{t\rightarrow \infty}r_i(t)=+\infty$.
For fixed $i\in V$, using Lemma \ref{limit behavior of theta_i}, we can choose $l>0$ large enough
such that, when $r_i>l$, the inner angle $\theta_i^{jk}$ is smaller than $\frac{\pi}{d_i}$, where
$d_i$ is the degree of the vertex $i$, which implies
$K_i=2\pi-\sum_{\triangle ijk\in F}\theta_i^{jk}>\pi>0$ and $R_i=\frac{K_i}{\tanh^2\frac{r_i}{2}}>0$.
Choose a time $t_0$ such that $r_i(t_0)>l$, this can be done as we have $\limsup_{t\rightarrow +\infty}r_i=+\infty$.
Set $a=\inf\{t<t_0|r_i(s)>l, \forall s\in [t, t_0]\}$, then $r_i(a)=l$.
Note that the hyperbolic Ricci flow (\ref{modifed hyperbolic Ricci flow}) could be written as
\begin{equation*}
\frac{dr_i}{dt}=\frac{1}{2}(\overline{R}_i-R_i)\sinh r_i.
\end{equation*}
So we have $r_i'(t)<0$ on $[a, t_0]$. Combining with $r_i(a)=l$, we have $r_i(t_0)<l$, which contradicts $r_i(t_0)>l$.

The arguments show that the solution $r(t)$ of the hyperbolic Ricci flow (\ref{modifed hyperbolic Ricci flow})
stays in a compact subset of $\mathds{R}^N_{>0}$
if the flow (\ref{modifed hyperbolic Ricci flow}) develops no removable singularities at finite time.
Thus $\widetilde{F}(u(t))$ is bounded and decreasing along
the extended hyperbolic flow (\ref{extended modified hyperbolic combinatorial Ricci flow}),
which implies $\lim_{t\rightarrow +\infty}\widetilde{F}(u(t))$ exist.
Then there exists $\xi_n\in (n, n+1)$ such that
$$\widetilde{F}(u(n+1))-\widetilde{F}(u(n+1))=-\sum_{i=1}^N(\overline{R}_i-\widetilde{R}_i)^2\tanh^2\frac{r_i}{2}|_{\xi_n}\rightarrow 0.$$
Since the solution $\{r(t)\}\subset\subset \mathds{R}^N_{>0}$, $r(\xi_n)$ have a convergent subsequence, still
denoted as $r(\xi_n)$, with $r(\xi_n)\rightarrow r^*$. Then we have $\widetilde{R}(r^*)=\overline{R}$.
By the uniqueness of $\widetilde{R}$ with respect to $r\in \mathds{R}^N_{>0}$ for the case of hyperbolic background geometry as
in Remark \ref{uniqueness in R^N}, we have $r^*=\overline{r}\in \Omega$.

The rest of the proof is an application of the Lyapunov Theorem to the hyperbolic combinatorial Ricci flow (\ref{modifed hyperbolic Ricci flow}).
Write the hyperbolic combinatorial Ricci flow (\ref{modifed hyperbolic Ricci flow}) in the form
\begin{equation}\label{hyperbolic flow in u}
\frac{du_i}{dt}=\overline{R}_i-R_i
\end{equation}
and set $\Gamma_i(u)=\overline{R}_i-R_i$. Then we have
$$D_u\Gamma|_{\overline{u}}=-\Sigma^{-1}L+\left(
                                            \begin{array}{ccc}
                                              \overline{R}_1 &   &   \\
                                                & \ddots &   \\
                                                &   & \overline{R}_N \\
                                            \end{array}
                                          \right),
$$
whose eigenvalues are all negative by Lemma \ref{positivity of L} and $\overline{R}\leq 0$.
Here $\Sigma=diag\{\tanh^2\frac{r_1}{2}, \cdots, \tanh^2\frac{r_N}{2}\}$.
This implies $\overline{u}$ is a local attractor of the hyperbolic combinatorial Ricci flow (\ref{hyperbolic flow in u}).
Note that $u(\xi_n)\rightarrow \overline{u}$, we get the solution $u(t)$ of (\ref{hyperbolic flow in u})
converges to $\overline{u}$ exponentially fast by Lyapunov Stability Theorem, which is equivalent to the solution $r(t)$
of (\ref{modifed hyperbolic Ricci flow}) converges to $\overline{r}$ exponentially fast.

  \item[(3)] The proof of the third part is a copy of the second part with Lemma \ref{limit behavior of theta_i}
replaced by Lemma \ref{limit behavior of extended theta}. We omit the details here.
\end{description}
\qed

%



\section{$\alpha$-curvature and $\alpha$-curvature flows}\label{section alpha curvature and alpha curvature flows}
Motivated by the generalization of the combinatorial curvature for Thurston's circle packing metrics to $\alpha$-curvatures in \cite{GX4},
we introduce the following definition of $\alpha$-curvature, which is a generalization of the combinatorial $R$-curvature in
Definition \ref{definition of comb curv with inversive dist}.

\begin{definition}
Given a closed triangulated surface $(M, \mathcal{T})$ with inversive distance $I\geq 0$
and a circle packing metric $r:V\rightarrow [0,+\infty)$, $\alpha\in \mathds{R}$ is a parameter,
the combinatorial $\alpha$-curvature $R_{\alpha}$ at the vertex $v_i$ is defined
to be
\begin{equation}\label{definition of alpha curv with Euclidean background}
R_{\alpha, i}=\frac{K_i}{s_i^\alpha}
\end{equation}
where $K_i$ is given by (\ref{classical Gauss curv}), $s_i(r)=r_i$ in the case of Euclidean background geometry
and $s_i(r)=\tanh\frac{r_i}{2}$ in the case of hyperbolic background geometry.
\end{definition}

\begin{remark}
In the case of $\alpha=0$, the $0$-curvature $R_{0}$ is just the classical combinatorial Gauss curvature $K$.
In the case of $\alpha=2$, the $2$-curvature $R_{2}$ is just the combinatorial $R$-curvature in Definition \ref{definition of comb curv with inversive dist}.
\end{remark}

\begin{definition}
We call a function $\overline{R}\in C(V)$ is $\alpha$-admissible if there is an inversive distance circle packing
metric $\overline{r}\in \Omega$ such that $R_\alpha(\overline{r})=\overline{R}$.
\end{definition}

Similar to Theorem \ref{theorem global rigidity}, we have the following global rigidity for $\alpha$-curvature $R_\alpha$.

\begin{theorem}\label{global rigidity for alpha curvature}
Given a weighted closed triangulated surface $(M,\mathcal{T}, I)$ with $I\geq0$. $\alpha\in \mathds{R}$ and $\overline{R}\in C(V)$
is a given function on $M$.
\begin{description}
  \item[(1)]  In the case of Euclidean background geometry,
  if $\alpha\overline{R}\equiv0$, there exists at most one inversive distance circle packing metric $\overline{r}\in \Omega$
  with combinatorial $\alpha$-curvature $\overline{R}$ up to scaling.
  For $\overline{R}\in C(V)$ with $\alpha\overline{R}\leq0$ and $\alpha\overline{R}\not\equiv0$,
  there exists at most one inversive distance circle packing metric $\overline{r}\in \Omega$ with combinatorial $\alpha$-curvature $\overline{R}$.
  \item[(2)] In the case of hyperbolic background geometry, if $\alpha\overline{R}\leq 0$,
  there exists at most one inversive distance circle packing metric $\overline{r}\in \Omega$
  with combinatorial $\alpha$-curvature $\overline{R}$.
\end{description}
\end{theorem}

Theorem \ref{global rigidity for alpha curvature} can be proved similarly to Theorem \ref{theorem global rigidity},
using the following Ricci potential function
\begin{equation*}
\begin{aligned}
F(u)=-\sum_{\Delta ijk\in F}F_{ijk}(u)+\int_{u_0}^u\sum_{i=1}^N(2\pi-\overline{R}_is_i^\alpha)du_i
\end{aligned}
\end{equation*}
defined on $\mathcal{U}$ and its convex extension
$$\widetilde{F}(u)=-\sum_{\triangle ijk\in F}\widetilde{F}_{ijk}(u)+\int_{u_0}^u\sum_{i=1}^N(2\pi-\overline{R}_is_i^\alpha)du_i,$$
where $u_i=\ln s_i$.
we omit the details here.

\begin{remark}
If we define $\widetilde{R}_{\alpha, i}=\frac{\widetilde{K}_i}{s_i^\alpha}$,
where $\widetilde{K}_i=2\pi-\sum_{\triangle ijk\in F}\widetilde{\theta}_i^{jk}$, we can also
get global rigidity results similar to that in Remark \ref{uniqueness in R^N}
under the same conditions in  Theorem \ref{global rigidity for alpha curvature}.
\end{remark}

\begin{remark}
The $\alpha$-curvature is obviously an extension of the classical combinatorial Gauss curvature $K$.
Bowers and Stephenson \cite{BS} once conjectured that the classical combinatorial Gauss curvature $K$ has global rigidity with respect to
the inversive distance circle packing metric $r$, which was proved by Guo \cite{Guo} and Luo \cite{L2} in the case of
Euclidean and hyperbolic background geometry. For the spherical background geometry, Ma and Schlenker \cite{MS} once
gave a counterexample to the conjecture, which is also available here.
Very recently, John C. Bowers and Philip L. Bowers presented an interesting elementary construction
of the Ma-Schlenker counterexample using only the inversive geometry of the 2-sphere \cite{BB}.
So we partially solves a generalized Bowers-Stephenson conjecture in Theorem \ref{global rigidity for alpha curvature}.
\end{remark}

We can also use the combinatorial Ricci flow to study the corresponding constant $\alpha$-curvature
problem and prescribing $\alpha$-curvature problem. As the proofs of these results are parallel to that
for the $R$-curvature in Section \ref{combinatorial Ricci flow}, we just list the results here and omit their proofs.

\begin{theorem}\label{constant alpha curvature problem}
Given a weighted closed triangulated surface $(M,\mathcal{T}, I)$ with $I\geq0$.
For the Euclidean background geometry, the normalized combinatorial $\alpha$-Ricci flow is defined to be
\begin{equation}\label{alpha Ricci flow}
\frac{dr_i}{dt}=(R_{\alpha, av}-R_{\alpha, i})r_i,
\end{equation}
where
\begin{equation*}
\begin{aligned}
R_{\alpha, av}=\left\{
                        \begin{array}{ll}
                          \frac{2\pi\chi(M)}{||r||_{\alpha}^\alpha}, & \hbox{$\alpha\neq 0$;} \\
                          \frac{2\pi\chi(M)}{N}, & \hbox{$\alpha=0$.}
                        \end{array}
                      \right.
\end{aligned}
\end{equation*}
\begin{description}
  \item[(1)] Along the normalized $\alpha$-Ricci flow (\ref{alpha Ricci flow}), the $\alpha$-curvature $R_\alpha$ evolves according to
$$\frac{dR_{\alpha, i}}{dt}=\Delta_{\alpha}R_{\alpha, i}+\alpha R_{\alpha, i}(R_{\alpha, i}-R_{\alpha, av}),$$
where the $\alpha$-Laplace operator is defined to be
\begin{equation*}
\Delta_{\alpha}f_i=\frac{1}{r_i^\alpha}\sum_{j=1}^N \left(-\frac{\partial K_i}{\partial u_j}\right) f_j
=\frac{1}{r_i^\alpha}\sum_{j\sim i} \left(-\frac{\partial K_i}{\partial u_j}\right) (f_j-f_i)
\end{equation*}
with $u_i=\ln r_i$ for $f\in C(V)$.
  \item[(2)] If the solution of the normalized $\alpha$-Ricci flow (\ref{alpha Ricci flow}) stays in $\Omega$ and
converges to an inversive distance circle packing metric $r^*\in \Omega$, then there exists a constant $\alpha$-curvature
metric in $\Omega$ and $r^*$ is such one.
  \item[(3)]Suppose $\alpha\chi(M)\leq 0$.
If there exists a Euclidean inversive distance circle packing metric $r^*\in \Omega$ with constant $\alpha$-curvature,
then the solution of the normalized $\alpha$-Ricci flow (\ref{alpha Ricci flow}) develops no essential singularities
in finite time and at time infinity; Furthermore, if the solution of (\ref{alpha Ricci flow}) develops no
removable singularities in finite time, then the solution of (\ref{alpha Ricci flow}) exists for all time,
converges exponentially fast to the constant $\alpha$-curvature metric $r^*\in \Omega$ and does not develop removable
singularities at time infinity.
  \item[(4)] If $\alpha\chi(M)\leq 0$ and there exists a Euclidean inversive distance
circle packing metric $r^*\in \Omega$ with constant $\alpha$-curvature, then the solution of the extended combinatorial
$\alpha$-Ricci flow
$$\frac{dr_i}{dt}=(R_{\alpha, av}-\widetilde{R}_{\alpha, i})r_i$$
exists for all time and converges exponentially fast to $r^*\in \Omega$.
\end{description}
\end{theorem}

We can also use the combinatorial $\alpha$-Ricci flow to study the prescribing $\alpha$-curvature problem and
get results similar to Theorem \ref{Euclidean prescribing curvature problem} by replacing $\overline{R}\leq 0$ with
$\alpha\overline{R}\leq 0$. As the results for the Euclidean and hyperbolic background geometry are similar, we state
them together in a unified form here.

\begin{theorem}\label{prescribling alpha curvature problem}
Given a closed triangulated surface $(M, \mathcal{T})$ with inversive distance $I\geq 0$.
Given a function $\overline{R}\in C(V)$, the modified combinatorial $\alpha$-Ricci flow is defined to be
\begin{equation}\label{modifed alpha Ricci flow}
\frac{ds_i}{dt}=(\overline{R}_i-R_{\alpha,i})s_i,
\end{equation}
where $s_i=r_i$ for the Euclidean background geometry and $s_i=\tanh\frac{r_i}{2}$ for the hyperbolic background geometry.
\begin{description}
  \item[(1)] If the solution of the modified hyperbolic $\alpha$-Ricci flow (\ref{modifed alpha Ricci flow})
stays in $\Omega$ and converges to $\overline{r}\in \Omega$, then we have $\overline{R}$ is $\alpha$-admissible and
$R_\alpha(\overline{r})=\overline{R}$.
  \item[(2)] Suppose $\overline{R}\in C(V)$ is $\alpha$-admissible with
$R_\alpha(\overline{r})=\overline{R}$ and $\alpha\overline{R}\leq 0$,
then the modified $\alpha$-Ricci flow (\ref{modifed alpha Ricci flow})
develops no essential singularities in finite time and at time infinity.
Furthermore, if the solution of the $\alpha$-Ricci flow (\ref{modifed alpha Ricci flow})
develops no removable singularities in finite time, then the solution of (\ref{modifed alpha Ricci flow})
exists for all time, converges exponentially fast to $\overline{r}$ and does not develop removable singularities at time infinity.
  \item[(3)] If
  $\overline{R}\in C(V)$ is $\alpha$-admissible with
  $R_{\alpha}(\overline{r})=\overline{R}$ for some $\overline{r}\in \Omega$ and $\alpha\overline{R}\leq 0$,
  then the solution of the extended $\alpha$-Ricci flow
  $$\frac{ds_i}{dt}=(\overline{R}_i-\widetilde{R}_{\alpha,i})s_i$$
  exists for all time and converges exponentially fast to $\overline{r}$.
\end{description}
\end{theorem}

\textbf{Acknowledgements}\\[8pt]
The authors thank the referees for their careful reading of the manuscript and for their nice suggestions.
The research of the first author is supported by National Natural Science Foundation of China (No.11501027)
and Fundamental Research Funds for the Central Universities (No. 2017JBM072).
The research of the second author is supported by National Natural Science Foundation of China under
grant no. 11301402.\\

\textbf{Note added in proof}\\[8pt]
Motivated by an observation of Zhou \cite{Z},
the second author \cite{X} recently proved the global rigidity of inversive distance circle packings in both Euclidean
and hyperbolic background geometry when the inversive distance $I>-1$ and satisfies
\begin{equation}\label{symmetric condition}
I_{ij}+I_{ik}I_{jk}\geq 0, I_{ik}+I_{ij}I_{jk}\geq 0, I_{jk}+I_{ij}I_{ik}\geq 0
\end{equation}
for any triangle $\triangle ijk$ in the triangulation of a surface.
The results obtained in this paper are all valid in the same form if the condition $I\geq 0$ is replaced by
$I>-1$ together with (\ref{symmetric condition}).

(Huabin Ge)
Department of Mathematics, Beijing Jiaotong University, Beijing 100044, P.R. China

E-mail: hbge@bjtu.edu.cn\\[2pt]

(Xu Xu) School of Mathematics and Statistics, Wuhan University, Wuhan 430072, P.R. China

E-mail: xuxu2@whu.edu.cn\\[2pt]

\end{document}